\documentclass[12pt]{amsart}
\usepackage{amscd, amsfonts, amssymb, amsthm}
\usepackage[all]{xypic}
\usepackage{graphics, lscape}
\usepackage{rotating}
\renewenvironment{enumerate}
  {\begin{list}
  {(\thesubsection.\theenumi)}
  {\usecounter{enumi}\setlength{\topsep}{8pt plus 2pt}
     \setlength{\parsep}{5.16663pt plus 1pt}
     \setlength{\itemsep}{5.16663pt plus 1pt}
     \setlength{\labelwidth}{2.5em}
     \setlength{\labelsep}{0.5em}
     \setlength{\itemindent}{0pt}
   }}
  {\end{list}}

\newcommand\CB{{\mathcal B}}

\newcommand\CN{{\mathcal N}} 

\newcommand\CP{{\mathcal P}}
\newcommand\CQ{{\mathcal Q}}

\newcommand\fg{{\mathfrak g}}

\newcommand\ft{{\mathfrak t}}

\newcommand\BBD{{\mathbb D}}
\newcommand\BBQ{{\mathbb Q}}

\newcommand\Ad{{\operatorname{Ad}}}

\newcommand\chk{^{\vee}}
\newcommand\CHom{{\mathcal{H}om}}
\newcommand\codim{{\operatorname{codim}}}
\newcommand\Ext{{\operatorname{Ext}}}
\newcommand\End{{\operatorname{End}}}

\newcommand\Hom{{\operatorname{Hom}}}

\newcommand\Lie{{\operatorname{Lie}}}

\newcommand\op{{^{\operatorname{op}}}}

\newcommand\inverse{^{-1}}
\renewcommand\th{{^{\text{th}}}}

\newcommand\Av{\operatorname{Av}}
\newcommand\bc{{\operatorname{bc}}}
\newcommand\bct{{\widetilde{\operatorname{bc}}}}
\newcommand\can{{\operatorname{can}}}
\newcommand\CNt{{\widetilde{\mathcal N}}}
\newcommand\CNtP{{{\widetilde{\mathcal N}}^\CP}}

\newcommand\Db{{D^b}}
\newcommand\dc{\operatorname{dc}}
\newcommand\Dbc{{D^b_c}}

\newcommand\dXPQcd{{\delta_{c,d}^{\CP,\CQ}}}

\newcommand\etabar{{\overline{\eta}}}

\newcommand\etarsP{\eta_{\rs}^\CP}

\newcommand\fgt{{\widetilde{\fg}}}
\newcommand\fgtP{{\widetilde{\fg}^{\CP}}}
\newcommand\fgtQ{{\widetilde{\fg}^{\CQ}}}
\newcommand\fgrs{\fg_{\rs}}
\newcommand\ftrs{\ft_{\reg}}
\newcommand\fgrst{{\widetilde {\fg}_{\rs}}}
\newcommand\fgrstP{{\fgrst^\CP}}
\newcommand\fgrstQ{\fgrst^\CQ}

\newcommand\HPQ{{H^{P,Q}}}

\newcommand\iM{i_{ {}_ M}}
\newcommand\iN{i_{ {}_ N}}
\newcommand\iP{i_{ {}_ P}}

\newcommand\IC{{\operatorname{IC}}}
\newcommand\ic{{\operatorname{ic}}}
\newcommand\id{{id}}

\newcommand\Jbar{\overline J}
\newcommand\jM{j_{ {}_ M}}
\newcommand\jN{j_{ {}_ N}}
\newcommand\jCN{j_{ {}_ {\CN}}}
\newcommand\jNt{j_{ {}_ {\CNt}}}
\newcommand\jNtP{j_{ {}_ {\CNtP}}}
\newcommand\jZ{j_{ {}_ Z}}
\newcommand\jX{j_{ {}_ X}}
\newcommand\jP{j_{ {}_ P}}

\newcommand\mubar{\overline{\mu}}

\newcommand\murs{{\mu_\rs}}

\newcommand\nat{{\operatorname{nat}}}
\newcommand\natt{\operatorname{nat}^\otimes}
\newcommand\nuM{\nu_{ {}_M}}
\newcommand\nuMr{\nu_{ {}_{ M_r}}}
\newcommand\nuP{\nu_{ {}_P}}

\newcommand\pr{{\operatorname{pr}}}

\newcommand\reg{{\operatorname{reg}}}
\newcommand\rs{{\operatorname{rs}}}

\newcommand\vd{\operatorname{vd}}
\newcommand\WP{{W_P}}
\newcommand\WQ{{W_Q}}
\newcommand\WPQ{{W^{P, Q}}}
\newcommand\XPQ{{X^{\CP,\CQ}}} 
\newcommand\xibar{\overline{\xi}}

\newcommand\xirsP{\xi_{\rs}^\CP}

\newcommand\XPQcd{{X^{\CP,\CQ}_{c,d}}}

\newcommand\YPQ{{Y^{\CP,\CQ}}} 
\newcommand\Ywbar{{\overline{Y_w}}}

\newcommand\ZPQ{{Z^{\CP,\CQ}}} 

\newcommand\Zwbar{{\overline{Z_w}}} 

\numberwithin{equation}{subsection}

\theoremstyle{plain}
\newtheorem{lemma}[equation]{Lemma}
\newtheorem{theorem}[equation]{Theorem}
\newtheorem{corollary}[equation]{Corollary}
\newtheorem{proposition}[equation]{Proposition}

\newtheorem*{lemma*}{Lemma}
\newtheorem*{lemmaA*}{Lemma A}
\newtheorem*{lemmaB*}{Lemma B}
\newtheorem*{corollary*}{Corollary}
\newtheorem*{corollaryB*}{Corollary B}
\newtheorem*{propositionA*}{Proposition A}
\newtheorem*{propositionC*}{Proposition C}
\newtheorem*{proposition*}{Proposition}
\newtheorem*{theorem*}{Theorem}

\thanks{Part of the research for this paper was carried out while both
  authors were staying at the Mathematisches Forschungsinstitut
  Oberwolfach supported by the ``Research in Pairs'' program.}

\thanks{Also part of this paper was written during visits of the first
  author to the University of Birmingham supported by an EPSRC grant
  to the second author.}


\begin{document}

\title[Homology of Steinberg Varieties] {Homology of generalized
Steinberg varieties\\ and\\ Weyl group invariants} 

\author[J.M. Douglass]{J. Matthew Douglass} \address{Department of
  Mathematics\\ University of North Texas\\ Denton TX, USA 76203}
\email{douglass@unt.edu} \urladdr{http://hilbert.math.unt.edu}

\author[G. R\"ohrle]{Gerhard R\"ohrle}
\address
{Fakult\"at f\"ur Mathematik, 
Ruhr-Universit\"at Bochum, 
 Universit\"atsstrasse 150, 
 D-44780 Bochum, Germany} \email{gerhard.roehrle@rub.de}
\urladdr{http://www.ruhr-uni-bochum.de/ffm/}


\maketitle
\allowdisplaybreaks

\begin{abstract}
  Let $G$ be a complex, connected, reductive algebraic group.  In this
  paper we show analogues of the computations by Borho and MacPherson
  of the invariants and anti-invariants of the cohomology of the
  Springer fibres of the cone of nilpotent elements, $\CN$, of
  $\Lie(G)$ for the Steinberg variety $Z$ of triples.
  
  Using a general specialization argument we show that for a parabolic
  subgroup $W_P \times W_Q$ of $W \times W$ the space of $W_P \times
  W_Q$-invariants and the space of $W_P \times W_Q$-anti-invariants of
  $H_{4n}(Z)$ are isomorphic to the top Borel-Moore homology groups of
  certain generalized Steinberg varieties introduced in
  \cite{douglassroehrle:geometry}.
  
  The rational group algebra of the Weyl group $W$ of $G$ is
  isomorphic to the opposite of the 
  top Borel-Moore homology $H_{4n}(Z)$ of $Z$, where
  $2n = \dim \CN$.  Suppose $W_P \times W_Q$ is a parabolic subgroup
  of $W \times W$.  We show that the space of $W_P \times
  W_Q$-invariants of $H_{4n}(Z)$ is $e_Q\BBQ We_P$, where $e_P$ is the
  idempotent in the group algebra of $W_P$ affording the trivial
  representation of $W_P$ and $e_Q$ is defined similarly.  We also
  show that the space of $W_P \times W_Q$-anti-invariants of
  $H_{4n}(Z)$ is $\epsilon_Q\BBQ W\epsilon_P$, where $\epsilon_P$ is
  the idempotent in the group algebra of $W_P$ affording the sign
  representation of $W_P$ and $\epsilon_Q$ is defined similarly.
\end{abstract}


\section{Introduction}

Suppose $G$ is a complex, reductive algebraic group and $\CB$ is the
variety of Borel subgroups of $G$. Then $\CB$ is a smooth, projective
variety. Let $T$ be a maximal torus in $G$ and choose a Borel
subgroup, $B$, of $G$ containing $T$. Let $W=N_G(T)/T$ be the Weyl
group of $(G,T)$. Then $W$ acts on $G/T$ on the right, the natural
projection $G/T\to G/B$ has the structure of a vector bundle, and the
varieties $G/B$ and $\CB$ are isomorphic. Thus, $W$ acts on the
singular cohomology with rational coefficients of $\CB$ via the
isomorphisms $H^\bullet(\CB)\cong H^\bullet(G/B)\cong H^\bullet(G/T)$.

Now suppose $P$ is a parabolic subgroup of $G$ containing $B$ and
$\CP$ is the variety of $G$-conjugates of $P$. Then $\CP$ is again a
smooth, projective variety and it is a classical result that
$H^\bullet (\CP)$ is isomorphic to the space of $W_P$-invariants in
$H^\bullet(\CB)$ where $W_P= N_P(T)/T$ is the Weyl group of $(P,T)$
(see \cite{hiller:geometry}).

Borho and MacPherson have generalized this result to fixed point
subvarieties of $\CB$ as follows. Let $\fg$ be the Lie algebra of $G$
and $\CN$ the cone of nilpotent elements in $\fg$. There is a
\emph{moment map}, $\mu_0\colon T^*\CB\to \CN$, where $T^*\CB$ is the
cotangent bundle of $\CB$.  For $x$ in $\CN$, set $\CB_x=\mu_0
\inverse(x)$. The variety $\CB_x$ may be identified with the variety
of all Borel subgroups of $G$, whose Lie algebra contains $x$. The
varieties $\CB_x$ vary from a point, when $x$ is regular, to $\CB$,
when $x=0$.  The moment map factors as $\mu_0=\eta_0 \circ \xi_0$
where $\xi_0\inverse(x)$ may be identified with the variety, $\CP_x$,
of all subgroups in $\CP$ whose Lie algebra contains $x$. There is
also a moment map, $\mu^{\CP}_0$ from the cotangent bundle of $\CP$ to
$\CN$, and $(\mu^{\CP}_0)\inverse (x)$ may be identified with the
variety of all subgroups in $\CP$ whose Lie algebras contain $x$ in
their nilradical. Set $\CP^0 _x= (\mu^{\CP}_0) \inverse(x)$.

Springer \cite{springer:trig} has defined an action of $W$ on
$H^\bullet(\CB_x)$ and Borho and MacPherson
\cite{borhomacpherson:partial} have shown that if $W$ acts on
$H^\bullet(\CB_x)$ by the tensor product of Springer's action with the
sign representation, then:
\begin{itemize}
\item[(1.1)] $H^\bullet(\CP_x)$ is isomorphic to the space of
  $W_P$-invariants in $H^\bullet(\CB_x)$.
\item[(1.2)] $H^\bullet(\CP_x^0)$ is isomorphic to the subspace of
  $H^\bullet(\CB_x)$ on which $W_P$ acts as the sign representation.
\end{itemize}

In a different direction, the \emph{Steinberg variety} of $G$ is the
fibred product $T^*\CB\times_\CN T^*\CB$ which may be identified with
the closed subvariety
\[
Z=\{\, (x, B', B'')\in \CN\times \CB\times \CB\mid x\in \Lie(B')\cap
\Lie(B'')\,\}
\]
of $\CN\times \CB\times \CB$. Kazhdan and Lusztig
\cite{kazhdanlusztig:topological} have defined an action of $W\times
W$, on $H_\bullet(Z)$, the rational, Borel-Moore homology of $Z$, and
they showed that the representation of $W\times W$ on the
top-dimensional homology group of $Z$, $H_{4n}(Z)$, where $n=\dim
\CB$, is equivalent to the two-sided regular representation of $W$.

Tanisaki \cite{tanisaki:twisted} and, more recently, Chriss and
Ginzburg \cite{chrissginzburg:representation} have strengthened the
connection between $H_\bullet(Z)$ and $W$ by defining a $\BBQ$-algebra
structure on $H_\bullet(Z)$ so that $H_i(Z) \cdot H_j(Z) \subseteq
H_{i+j-4n}(Z)$ and $H_{4n}(Z)\op$ is isomorphic to the group algebra
$\BBQ W$.

In this paper we prove analogs of (1.1) and (1.2) for the Steinberg
variety.

Suppose $Q$ is a parabolic subgroup of $G$ containing $B$ (a special
case is when $Q=P$), $W_Q$ is the Weyl group of $(Q,T)$, and $\CQ$ is
the conjugacy class of parabolic subgroups that contains $Q$. In
\cite{douglassroehrle:geometry} we defined generalized Steinberg
varieties
\[
\XPQ= \{\, (x, P', Q') \in \CN \times \CP \times \CQ \mid x\in
\Lie(P') \cap \Lie(Q') \,\}
\]
and
\[
\YPQ= \{(x, P', Q') \in \CN \times \CP \times \CQ \mid x\in
\Lie(U_{P'}) \cap \Lie(U_{Q'}) \,\}
\]
where $U_{P'}$ and $U_{Q'}$ are the unipotent radicals of $P'$ and
$Q'$ respectively. It was shown in \cite{douglassroehrle:geometry}
that $\XPQ$ is purely $2n$-dimensional and $\YPQ$ is purely
$(2n-f)$-dimensional where $f= \dim P/B +\dim Q/B$.

The first analogs of (1.1) and (1.2) are:
\begin{itemize}
\item[(1.1$'$)] $H_{4n}(\XPQ)$ is isomorphic to the space of
  $W_P\times W_Q$-invariants in $H_{4n}(Z)$.
\item[(1.2$'$)] $H_{4n-2f}(\YPQ)$ is isomorphic to the subspace of
  $H_{4n}(Z)$ on which $W_P\times W_Q$ acts as the sign
  representation.
\end{itemize}
We prove both of these statements in this paper.

More generally we consider the following statements:
\begin{itemize}
\item[(1.1$''$)] $H_\bullet(\XPQ)$ is isomorphic to the space of
  $W_P\times W_Q$-invariants in $H_\bullet(Z)$.
\item[(1.2$''$)] $H_\bullet(\YPQ)$ is isomorphic to the subspace of
  $H_\bullet(Z)$ on which $W_P\times W_Q$ acts as the sign
  representation.
\end{itemize}

In \S3 we prove a general specialization result, in the spirit of
\cite{borhomacpherson:partial}, which has (1.1$''$) as a special case.
Obviously (1.1$'$) follows immediately from (1.1$''$).  It seems
likely that (1.2$''$) is true, but our proof of (1.2$'$) uses
dimension computations from \cite{douglassroehrle:geometry} that are
not available for $H_i(\YPQ)$ for $i<4n-2f$.

In \S4 we prove a general equivariance result in the spirit of
\cite{chrissginzburg:representation}. A special case of this result is
that there is a $W\times W$-equivariant isomorphism
\begin{equation*}
  \xymatrix{\Ext_{\CN}^{4n-\bullet} \left( R(\mu_0)_! \BBQ_{T^*\CB},
      R(\mu_0)_! \BBQ_{T^*\CB} \right) \ar[r]^-{\simeq} &
    H_{\bullet}(Z)\op.} 
\end{equation*}
Borho and MacPherson \cite{borhomacpherson:weyl} have shown that the
$\BBQ$-algebras $\BBQ W$ and $\End_{\CN}( R(\mu_0)_!  \BBQ_{T^*\CB})$
are isomorphic and Chriss and Ginzburg \cite[\S8.6]
{chrissginzburg:representation} have shown that 
\[
\Ext_{\CN}^{4n-\bullet} \left( R(\mu_0)_! \BBQ_{T^*\CB}, R(\mu_0)_!
  \BBQ_{T^*\CB} \right) \cong H_{\bullet}(Z)\op.
\]
Thus, taking $\bullet =4n$ we get $W\times W$-equivariant,
$\BBQ$-algebra isomorphisms
\begin{equation*}
\xymatrix{\BBQ W\ar[r]^-{\simeq} &\End_{\CN}( R(\mu_0)_! \BBQ_{T^*\CB})
  \ar[r]^-{\simeq}& H_{4n}(Z)\op.}
\end{equation*}
where $W\times W$ acts on $\BBQ W$ by $(w,w')\cdot v= w' v w\inverse$
for $w$ and $w'$ in $W$ and $v$ in $\BBQ W$.

Using the isomorphism between $\BBQ W$ and $H_{4n} (Z)\op$ we may
formulate (1.1$'$) and (1.2$'$) in terms of the group algebra of $W$:
\begin{itemize}
\item[(1.1$'''$)] If $e_P$ is the primitive idempotent in $\BBQ\WP$
  corresponding to the trivial representation of $\WP$ and $e_Q$ is
  defined similarly, then $H_{4n}(\XPQ)$ is isomorphic to the subspace
  $e_Q\BBQ We_P$ of $\BBQ W$.
\item[(1.2$'''$)] If $\epsilon_P$ is the primitive idempotent in
  $\BBQ\WP$ corresponding to the sign representation of $\WP$ and
  $\epsilon_Q$ is defined similarly, then $H_{4n-2f} (\YPQ)$ is
  isomorphic to the subspace $\epsilon_Q\BBQ W \epsilon_P$ of $\BBQ
  W$.
\end{itemize}

In \cite{douglassroehrle:geometry} we defined generalized Steinberg
varieties $\XPQcd$. Statements (1.1$'''$) and (1.2$'''$) together with
computations in some special cases suggest that the Borel-Moore
homology of a general $\XPQcd$ is given as follows.

A generalized Steinberg variety, $\XPQcd$, depends on a pair of
nilpotent adjoint orbits in $\Lie(P/U_P)$ and $\Lie(Q/U_Q)$
respectively. We will not recall the precise definition here but
instead refer the interested reader to
\cite{douglassroehrle:geometry}.  In turn, these nilpotent orbits
determine irreducible representations of $W_P$ and $W_Q$, say $\rho_c$
and $\rho_d$ respectively, corresponding to the trivial representation
of the component groups of the orbits via the Springer correspondence
as defined in \cite{borhomacpherson:weyl}. Let $e_c$ and $e_d$ denote
primitive idempotents in $\BBQ\WP$ and $\BBQ\WQ$ affording $\rho_c$
and $\rho_d$ respectively. In \cite[Corollary 2.6]
{douglassroehrle:geometry} we have given a sharp upper bound,
$\dXPQcd$, for the dimension of $\XPQcd$. We conjecture that
\begin{itemize}
\item $H_{2\dXPQcd}( \XPQcd)$ is isomorphic to $e_d\BBQ W e_c$.
\end{itemize}
More generally, we conjecture that
\begin{itemize}
\item $H_\bullet(\XPQcd)$ is isomorphic to $e_d H_\bullet(Z) e_c$
  where we consider $e_c$ and $e_d$ in $H_\bullet(Z)$ via the
  isomorphism $\BBQ W\cong H_{4n}(Z)\op$.
\end{itemize}

In much of this paper (\S2 -- \S4 and the Appendix) we are concerned
with general sheaf theory. Most of our conclusions about the
Borel-Moore homology of generalized Steinberg varieties are
straightforward applications of more general results. The main
theorems, which are described briefly below, are the specialization
results, Theorem \ref{thm3.5} and Corollary \ref{cor3.5}, and the
equivariance results discussed in \S\ref{4.1}.  We hope these general
results will have applications outside the realm of generalized
Steinberg varieties.

Our computation of the Borel-Moore homology of $\XPQ$ and $\YPQ$ is
given in \S5. Although the results depend on facts proved in \S3 and
\S4, this section may be read independently of the other sections.

The rest of this paper is organized as follows.

In \S2 we fix notation and collect some sheaf-theoretic results that
are used in subsequent sections for which we could not find a suitable
reference.

In \S3 we give an axiomatic approach to a specialization result which
allows us to identify a direct image map in Borel-Moore homology with
the averaging map for a group action. The basic idea goes back to
Lusztig \cite{lusztig:green} and Borho-MacPherson
\cite{borhomacpherson:partial}. A result which is similar in spirit,
but which is in a sense dual to our result, and does not apply to
Borel-Moore homology, has been used by Spaltenstein in
\cite{spaltenstein:reflection}. Statement (1.1$''$) is a
straightforward consequence of the main result in this section,
Theorem \ref{thm3.5}.

In \S4 we continue the axiomatic approach from \S3 and prove an
equivariance result for two-sided group actions that is key for our
application to generalized Steinberg varieties. The crucial result is
Theorem \ref{thm4.4} which when applied to the Steinberg variety
implies that there is a $W\times W$-equivariant isomorphism between
$\Ext_\CN^{4n-\bullet}(R(\mu_0)_! \BBQ_{T^*\CB},R(\mu_0)_! \BBQ_{T^*\CB})$ 
and $H_\bullet(Z)$.
This result is similar in spirit to the results in \cite[\S8.6]
{chrissginzburg:representation}.

In \S5 we specialize the results in the previous sections to the case
of generalized Steinberg varieties and prove (1.1$''$), (1.2$'$),
(1.1$'''$), and (1.2$'''$).

In the Appendix, we prove two results about the natural transformation
$\xi^*\to \xi^![2l]$ for a morphism $\xi\colon X\to Y$, where $l= \dim
Y- \dim X$. These results are needed in the proof of Theorem
\ref{thm4.4}.

For simplicity, in this paper we have chosen to work with complex
algebraic groups and Borel-Moore homology, but our arguments are
essentially categorical and make sense in the setting of algebraic
groups over arbitrary algebraically closed fields and $l$-adic
cohomology.

\section{Preliminaries}

\subsection{}\label{2.1}
First, we fix some assumptions and notation that will be used
throughout the rest of this paper. The reader is urged to skim this
section quickly to become familiar with the notation and refer back to
the results used in the sequel when necessary. The main references for
sheaf-theoretic notation and results used in this paper are the
article \cite{borel:intersection} by Borel (with the collaboration of
N. Spaltenstein) and the book \cite{kashiwarashapira:sheaves} by
Kashiwara and Shapira.

The topological spaces we consider are complex algebraic varieties
endowed with their Euclidean topologies, although many arguments apply
as well to pseudomanifolds as defined in \cite[\S1.1]
{goreskymacpherson:intersectionII}.

The ``dimension'' of a space always means its dimension as a complex
algebraic variety.

If $X$ is a variety, then $D(X)$ denotes the derived category of the
category of sheaves of $\BBQ$-vector spaces on $X$, $\Db(X)$ denotes
the full subcategory of $D(X)$ consisting of complexes with bounded
cohomology, and $\Dbc(X)$ denotes the full subcategory of $\Db(X)$
consisting of complexes with constructible cohomology.

For complexes $A$ and $B$ in $D(X)$, $\Ext^{j}( A,B)$ is
defined to be $H^j( R\Hom(A, B))$ and it is shown in
\cite[\S5.17]{borel:intersection} that $\Ext^{j}( A, B) =
\Hom_{D(X)}(A, B[j])$. Define
\[
\Ext_X^{j}( A, B) =\Hom_{D(X)}(A, B[j]).
\]

Since $\Dbc(X)$ is a full subcategory of $D(X)$, if $A$ and $B$ are
complexes in $\Dbc(X)$, then $\Hom_{\Dbc(X)}(A, B)= \Hom_{D(X)}(A,
B)$. To simplify the notation, we denote both of these spaces by
$\Hom_{X}(A, B)$. Also, we denote the complex $A\overset{L}{\otimes}
B$ simply by $A\otimes B$.

The constant sheaf on $X$, considered as a complex concentrated in
degree $0$, is denoted by $\BBQ_X$ and the dualizing complex of $X$ is
denoted by $\BBD_X$.

If $A$ is a complex of sheaves of $\BBQ$-vector spaces on $X$, then
$A\chk= R\CHom( A, \BBD_X)$ denotes the Verdier dual of $A$. There is
a canonical isomorphism between $\BBD_X$ and $\BBQ_X\chk$ we denote by
$\dc_X$, so
\[
\xymatrix{\dc_X\colon \BBD_X \ar[r]^-{\simeq}& \BBQ_X\chk}.
\]

If $f\colon A\to B$ is a morphism in $D(X)$, and $C$ is a complex in
$D(X)$, then $f$ induces natural morphisms in $D(X)$,
\[
\xymatrix{f^\sharp \colon R\CHom(B,C) \ar[r]& R\CHom(A,C)}
\quad\text{and}\quad \xymatrix{f_\sharp \colon R\CHom(C,A) \ar[r]&
  R\CHom(C,B).}
\]
In the special case when $C=\BBD_X$, we have $ R\CHom(A,C)= A\chk$ and
$R\CHom(B,C)= B\chk$. We usually write $f\chk$ instead of $f^\sharp$
in this case, so $f\chk \colon B\chk\to A\chk$ is the Verdier dual of
$f$.

Similarly, $f$ induces natural linear transformations
\[
\xymatrix{f^\sharp \colon \Ext^{\bullet}_X(B,C) \ar[r]&
  \Ext^{\bullet}_X(A,C)} \quad\text{and}\quad \xymatrix{f_\sharp
  \colon \Ext^{\bullet}_X(C,A) \ar[r]& \Ext^{\bullet}_X(C,B).}
\]

The $j\th$ Borel-Moore homology group of a locally compact, Hausdorff
topological space, $X$, has several equivalent definitions (see
\cite[\S2.6] {chrissginzburg:representation}). In this paper we use
the canonical isomorphisms,
\[
H^{-j}(X, \BBD_X) \cong H^{-j}(X, R\CHom( \BBQ_X, \BBD_X))\cong
\Ext^{-j}_X( \BBQ_X, \BBD_X)
\]
where $H^{-j}(X, \BBD_X)$ is the hypercohomology of $X$ with
coefficients in $\BBD_X$ and we \emph{define} the $j\th$ Borel-Moore
homology group of $X$ by
\[
H_j(X)= \Ext^{-j}_X( \BBQ_X, \BBD_X).
\]

\subsection{}\label{2.2}

Now suppose that $\xi\colon X\to Y$ is a morphism of varieties. Then
$\xi$ determines natural isomorphisms
\[
\xymatrix{\phi_\xi\colon R\CHom( R\xi_!A, B)\ar[r]^-{\simeq}& R\xi_*
  R\CHom( A, \xi^!B)}
\]
and
\[
\xymatrix{\nat_\xi\colon \xi^! R\CHom(B,C) \ar[r]^-{\simeq}& R\CHom(
  \xi^*B, \xi^! C)}
\]
for $A$ in $D(X)$ and $B$ and $C$ in $D(Y)$.

There are canonical isomorphisms, 
\[
\xymatrix{\alpha_\xi\colon \xi^* \BBQ_Y \ar[r]^-{\simeq}& \BBQ_X}
\qquad \text{and} \qquad \xymatrix{\beta_\xi \colon \BBD_X
  \ar[r]^-{\simeq}& \xi^!  \BBD_Y}
\]
in the category of sheaves on $X$ and $\Dbc(X)$ respectively. It is
straightforward to check that $\alpha_\xi$ and $\beta_\xi$ have the
following properties:
\begin{enumerate}
\item The maps $\beta_\xi \colon \BBD_X\to \xi^!  \BBD_Y$ and
  $(\beta_\xi)_\sharp \colon R\CHom(\xi^* \BBQ_Y, \BBD_X)
  \longrightarrow R\CHom( \xi^*\BBQ_Y, \xi^! \BBD_Y)$ are related by
  $(\beta_\xi)_\sharp \circ \alpha_\xi\chk \circ \dc_X= \nat_\xi \circ
  \xi^!(\dc_Y) \circ \beta_\xi$ where $\dc_X$ and $\dc_Y$ are as in
  \S\ref{2.1}.
  \label{2.2.1} 
\item If $\eta\colon Y\to Z$ is another morphism of varieties, then
  $\alpha_{\eta\xi}= \alpha_\xi\circ \xi^*( \alpha_\eta)$ and
  $\beta_{\eta\xi}= \xi^!(\beta_\eta)\circ \beta_\xi$.  \label{2.2.2}
\end{enumerate}

\subsection{}\label{2.3}
Let $\delta\colon X\to X\times X$ be the diagonal embedding and let
$p$ and $q$ denote the projections of $X\times X$ on the first and
second factors respectively.  In \cite[Theorem
10.25]{borel:intersection} it is shown that there is a natural
isomorphism,
\[
\xymatrix{\lambda \colon A\chk \boxtimes B \ar[r]^-{\simeq}&
  R\CHom(p^*A, q^!B)}
\]
in $\Dbc(X\times X)$. It follows that $\nat_\delta\circ
\delta^!(\lambda)$ is a natural isomorphism between $\delta^!(A\chk
\boxtimes B)$ and $R\CHom(A, B)$.

\begin{proposition}
  Suppose $A$ and $B$ are in $\Dbc(X)$, $u\colon A\to A$ is an
  endomorphism of $A$, and $v\colon B\to B$ is an endomorphism of $B$.
  Then the diagram
  \[
  \xymatrix{ A\chk\boxtimes B \ar[rrr]^{u\chk\boxtimes v}
    \ar[d]_\lambda &&&
    A\chk\boxtimes B \ar[d]^\lambda \\
    R\CHom( p^*A, q^!B) \ar[rrr]_{ (p^*u)^\sharp\circ (q^!v)_\sharp}
    &&& R\CHom( p^*A, q^!B) }
  \]
  commutes.
\end{proposition}

\begin{proof}
  By definition $A\chk \boxtimes B= p^*R\CHom(A,\BBD_X) \otimes q^*B$
  and $u\chk \boxtimes v=p^*(u^\sharp)\otimes q^*v$.
  
  In the special case when $A$ is the constant sheaf, the isomorphism
  $\lambda$ may be identified with a natural isomorphism $\lambda'
  \colon p^*\BBD_X\otimes q^*B\to q^!B$ as in \cite[\P10.24]
  {borel:intersection}.  Then for an arbitrary $A$, the isomorphism
  $\lambda$ is defined as the composition $\lambda' _\sharp \circ
  h_2\circ h_1$, where $h_1$ and $h_2$ are the natural maps
  \[
  \xymatrix{h_1\colon p^*R\CHom(A,\BBD_X) \otimes q^*B \ar[r]&
    R\CHom(p^*A, p^*\BBD_X) \otimes q^*B}
  \]
  and
  \[
  \xymatrix{h_2\colon R\CHom(p^*A, p^*\BBD_X) \otimes q^*B\ar[r]&
    R\CHom(p^*A, p^*\BBD_X \otimes q^*B).}
  \]
  
  It is straightforward to check that 
  \[
  h_1\circ \left( p^*(u^\sharp)\otimes q^*v \right) =\left(
    (p^*u)^\sharp)\otimes q^*v \right) \circ h_1
  \]
  and
  \[
  h_2\circ \left( (p^*u)^\sharp)\otimes q^*v \right) = \left(
    (p^*u)^\sharp \circ (id\otimes q^*v)_\sharp \right) \circ h_2.
  \]
  Moreover, it follows from the naturality of $\lambda'$ that
  \[
  \lambda' _\sharp \circ \left( (p^*u)^\sharp \circ (id\otimes
    q^*v)_\sharp \right) = \left( (p^*u)^\sharp \circ (q^!v)_\sharp
  \right) \circ \lambda' _\sharp.
  \]
  
  Therefore $\lambda \circ \left( u\chk\boxtimes v\right) = \left(
    (p^*u)^\sharp \circ (q^!v)_\sharp \right) \circ \lambda$, as
  desired.
\end{proof}

\begin{corollary}\label{cor2.3}
  With the preceding notation, the diagram
  \[
  \xymatrix{\delta^!\left(A\chk\boxtimes B\right) \ar[rr]^{\delta^!(
      u\chk\boxtimes v)} \ar[d]_{\nat_\delta \circ
      \delta^!(\lambda)}&& \delta^! \left( A\chk\boxtimes B\right)
    \ar[d]^{\nat_\delta \circ \delta^!(\lambda)} \\
    R\CHom( A,B) \ar[rr]_{u^\sharp \circ v_\sharp}&& R\CHom(A,B) }
  \]
 commutes.
\end{corollary}

\begin{proof}
  We have just seen that $\lambda \circ \left( u\chk\boxtimes v\right)
  = \left( (p^*u)^\sharp \circ (q^!v)_\sharp \right) \circ \lambda$,
  so
  \[
  \delta^!(\lambda) \circ \delta^!  (u\chk\boxtimes v) = \delta^!
  \left( (p^*u)^\sharp \circ (q^!v)_\sharp \right) \circ
  \delta^!(\lambda).
  \]
  It is straightforward to check that
  \[
  \nat_\delta\circ \delta^! \left( (p^*u)^\sharp \circ (q^!v)_\sharp
  \right)= \left ((\delta^*p^*u)^\sharp \circ (\delta^!q^!v)_\sharp
  \right) \circ \nat_\delta =\left ( u^\sharp \circ v_\sharp \right)
  \circ \nat_\delta
  \]
  so
  \[
  \nat_\delta \circ \delta^!(\lambda) \circ \delta^!  (u\chk\boxtimes
  v) = \nat_\delta \circ \delta^!  \left((p^*u)^\sharp \circ
    (q^!v)_\sharp \right) \circ \delta^!(\lambda)= \left (u^\sharp
    \circ v_\sharp \right) \circ \nat_\delta\circ \delta^!(\lambda).
  \]
  This proves the corollary.
\end{proof}

It is shown in \cite[\S2.6]{kashiwarashapira:sheaves} that for $A$,
$B$, and $C$ in $D(X)$ there is a natural isomorphism
$\Hom_{X}(C\otimes A, B) \cong \Hom_{X}(C,R\CHom(A, B))$. It follows
that there is an isomorphism of graded vector spaces
$\Ext_X^\bullet(C\otimes A, B) \cong \Ext_X^\bullet(C,R\CHom(A, B))$.
Taking $C=\BBQ_X$ and using the canonical isomorphism $\BBQ_X\otimes
A\cong A$ we get a natural isomorphism of graded vector spaces
\[
\xymatrix{\can \colon\Ext_X^\bullet(A, B) \ar[r]^-{\simeq}&
  \Ext_X^\bullet(\BBQ_X,R\CHom(A, B)).}
\]

The next proposition follows from the naturality of $\can$.

\begin{proposition} \label{prop2.3}
  Suppose $A$ and $B$ are in $D(X)$, $u\colon A\to A$ is an
  endomorphism of $A$, and $v\colon B\to B$ is an endomorphism of $B$.
  Then the diagram
  \[
  \xymatrix{ \Ext_X^\bullet(A, B) \ar[rr]^-{\can} \ar[d]_{u^\sharp \circ
      v_\sharp} && \Ext_X^\bullet(\BBQ_X,R\CHom(A, B))
    \ar[d]^{(u^\sharp \circ v_\sharp)_\sharp}\\
    \Ext_X^\bullet(A, B) \ar[rr]_-{\can}&&
    \Ext_X^\bullet(\BBQ_X,R\CHom(A, B))}
  \]
  commutes.
\end{proposition}

\subsection{}\label{2.4}
As in \S\ref{2.2}, $\xi\colon X\to Y$ is a morphism of varieties. The
functors $\xi^*$ and $R\xi_*$ form an adjoint pair. We denote by
\[
\xymatrix{\Psi_\xi\colon \Hom_{X} (\xi^*B,A) \ar[r]^-{\simeq}
  &\Hom_{Y}(B,R\xi_*A)}
\]
the adjunction mapping for $A$ in $D(X)$ and $B$ in $D(Y)$ and by
$\chi^\xi$ the unit of the adjunction. Although $\chi^\xi$ is a
natural transformation, $\chi^\xi_B\colon B \to R\xi_*\xi^* B$, in
order to simplify the notation we omit the subscript and just write
$\chi^\xi$ instead of $\chi^\xi_B$. The appropriate subscript is
always uniquely determined by the context and so this should cause no
confusion.

Similarly, the functors $R\xi_!$ and $\xi^!$ form an adjoint pair. We
denote by 
\[
\xymatrix{\Phi_\xi \colon \Hom_{Y}(R\xi_!A,B)\ar[r]^-{\simeq}
  &\Hom_{X}(A, \xi^!B)}
\]
the adjunction mapping and by $\epsilon^\xi$ the counit of the
adjunction.

We need the following identities for morphisms $f\colon R\xi_!A\to B$
and $k\colon B\to B'$ in $D(Y)$ and $g\colon A\to \xi^!B$ and $h\colon
A'\to A$ in $D(X)$ (see \cite[IV.1]{maclane:categories}):
\begin{gather}
  \epsilon^\xi=\Phi_\xi\inverse(\id) \qquad \Phi_\xi\inverse(g)=
  \epsilon^\xi\circ R\xi_!(g)\label{eq:ad1}\\
  \Phi_\xi(f\circ R\xi_!(h))= \Phi_\xi(f)\circ h \qquad \Phi_\xi(
  k\circ f)= \xi^!(k)\circ \Phi_\xi(f) \label{eq:ad2}\\
  \Phi_\xi\inverse(g \circ h)= \Phi_\xi\inverse(g)\circ R\xi_!(h)
  \qquad \Phi_\xi\inverse( \xi^!(k) \circ g)=k\circ \Phi_\xi
  \inverse(g)
  \label{eq:ad3}
\end{gather}

Verdier duality defines contravariant automorphisms of the
subcategories $\Dbc(X)$ and $\Dbc(Y)$ of $D(X)$ and $D(Y)$
respectively. In these subcategories we can use standard identities
for Verdier duality in \cite[\S10]{borel:intersection} to express
$\Phi_\xi$ and $\epsilon^\xi$ in terms of $\Psi_\xi$ and $\chi^\xi$ as
follows.

Suppose $A$ is in $\Dbc(X)$, $B$ is in $\Dbc(Y)$, and $f$ is in
$\Hom_{Y} (R\xi_!  A,B)$. Then $\Psi_\xi\inverse (f\chk)\chk$ is in
$\Hom_{X} ( A, \xi^!B)$. Clearly, $f\mapsto \Psi_\xi\inverse
(f\chk)\chk$ is natural in $A$ and $B$ and so we may \emph{define}
$\Phi_\xi$ by $\Phi_\xi(f)= \Psi_\xi\inverse (f\chk)\chk$.

Similarly, taking the Verdier dual of $\chi^\xi_B \colon B \to
R\xi_*\xi^* B$ we get $\left( \chi^\xi_B \right)\chk \colon R\xi_!
\xi^! B\chk \to B\chk$ and we conclude that $\left( \chi^\xi_B
\right)\chk = \epsilon^\xi_{B\chk}$.

\subsection{}\label{2.5}
Next, consider a cartesian square
\begin{equation}
  \label{eq:cart}
  \xymatrix{
    X' \ar[r]^{i} \ar[d]_{\eta}& X \ar[d]^{\xi}\\
    Y' \ar[r]^{j}& Y}  
\end{equation}
where $\xi$ and $\eta$ are proper morphisms. Then
$\Psi_j\inverse(R\xi_*(\chi^i)) \colon j^* R\xi_*\to R\eta_* i^*$ is a
natural equivalence of functors from $D(X)$ to $D(Y)$. Restricting to
$\Dbc(X)$ and $\Dbc(Y)$ and taking the Verdier dual we conclude that
$\Psi_j\inverse(R\xi_*(\chi^i)) \chk \colon R\eta_!  i^!\to j^!
R\xi_!$ is a natural equivalence. It follows from the discussion in
\S\ref{2.4} above that
\[
\Psi_j\inverse(R\xi_*(\chi^i)) \chk = \Phi_j\left(
  (R\xi_*(\chi^i))\chk \right) = \Phi_j\left( R\xi_!((\chi^i)\chk)
\right) = \Phi_j\left( (R\xi_!(\epsilon^i) \right).
\]
Define 
\[
\xymatrix{\bc_{\eta, i}\colon j^!  \circ R\xi_! \ar[r]& R\eta_!  \circ
  i^!} \quad \text{by} \quad \bc_{\eta, i} = \Phi_j \left( R\xi_!(
  \epsilon^i) \right)\inverse.
\]
Then $\bc_{\eta, i}$ is a natural equivalence and
$\bc_{\eta,i}\inverse = \Phi_j \left( R\xi_!( \epsilon^i) \right) $.

\begin{lemma}\label{lem2.5}
  Suppose that in diagram (\ref{eq:cart}) the maps $i$ and $j$ are
  open embeddings. Then, for $A$ in $\Dbc(X)$ and $B$ in $\Dbc(Y)$,
  the diagram
  \[
  \xymatrix{ {j^!R\xi_! R\CHom( A, \xi^! B)} \ar[rr]^{j^!(\phi_\xi
      \inverse)} \ar[d]_{\bc} && j^! R\CHom( R\xi_!A, B)
    \ar[r]^{\nat_j} & R\CHom(j^* R\xi_! A, j^! B)
    \ar[d]^{(\bc\inverse)^\sharp} \\
    R\eta_! i^! R\CHom( A, \xi^!B) \ar[rr]_{ R\eta_!(\nat_i)} &&
    R\eta_!  R\CHom( i^*A, i^! \xi^! B) \ar[r]_{\phi_\eta\inverse} &
    R\CHom( R\eta_! i^!A, j^!B) }
  \]
  commutes in $\Dbc(Y')$, where $\bc= \bc_{\eta, i}$.
\end{lemma}

\begin{proof}
  Since $i$ and $j$ are open embeddings, we have $i^!=i^*$ and $j^!=
  j^*$, so the statement of the lemma makes sense and is easily proved
  for sheaves on $X$ and $Y$. The result then follows using standard
  arguments for derived functors.
\end{proof}

\subsection{}\label{2.6}

If $U$ is a smooth, open, dense subvariety of $X$, and $L$ is a local
system on $U$, then we denote the intersection complex, as in
\cite{borel:intersection}, middle perversity, by $\IC(X, L)$. It is a
complex of sheaves in $\Dbc(X)$.  It is shown in \cite[Theorem
3.5]{goreskymacpherson:intersectionII} that $\IC$ defines a fully
faithful functor from the category of local systems on $U$ to
$\Dbc(X)$.

Notice that if we start with a complex, $A$, on an open, dense
subvariety of $X$ with $H^p(A)=0$ for $p\ne0$, then we may construct a
complex $\IC(X, A)$ as in \cite[\S2.2]{borel:intersection} starting
with $A$.  The complexes $\IC(X,A)$ and $\IC(X, H^0(A))$ are
isomorphic in $\Dbc(X)$.


\section{Specialization}

\subsection{}\label{3.1}
In this section we axiomatize a specialization argument that allows us
to compute invariants in Borel-Moore homology. There are various
schemes that allow one to use generic information to prove
(co)-homological results about special fibres, or more generally
closed subvarieties (see \cite{fultonmacpherson:categorical},
\cite{chrissginzburg:representation}, \cite{rossmann:picard}). Our
approach, which is based on an idea of Lusztig in \cite{lusztig:green}
that was generalized by Borho and MacPherson
\cite{borhomacpherson:partial}, is to use intersection complexes of
local systems on open, dense subvarieties of a variety, $N$, to obtain
information about the Borel-Moore homology groups of a closed
subvariety, $N_0$, of $N$.

We start with what we call the ``basic commutative diagram'' of
morphisms of complex, algebraic varieties consisting of cartesian
squares:
\begin{equation}
  \label{eq:bd}
  \xymatrix{ 
    M_0 \ar[rr]^{\eta_0} \ar[d]_{\jM}&& P_0 \ar[rr]^{\xi_0}
    \ar[d]_{\jP}&& N_0 \ar[d]_{\jN} \\
    M \ar[rr]^{\eta}&& P \ar[rr]^{\xi}&& N\\
    M_r \ar[rr]^{\eta_r} \ar[u]^{\iM}&& P_r \ar[rr]^{\xi_r}
    \ar[u]^{\iP}&& N_r \ar[u]^{\iN}}
\end{equation}
Define 
\[
\mu= \xi\eta, \quad \mu_r= \xi_r\eta_r,\quad \text{and} \quad \mu_0=
\xi_0\eta_0. 
\]

We assume that this basic commutative diagram has the following
properties:
\begin{itemize}
\item[D1] The varieties $M$, $P$, and $N$ are purely $d$-dimensional.
\item[D2] The varieties $M$, $P$, and $N_r$ are rational homology manifolds.
\item[D3] The morphisms $\xi$ and $\mu$ are surjective, proper
  morphisms that are \emph{small} (see
  \cite[\S6.2]{goreskymacpherson:intersectionII}) in the sense that
  for all $r>0$,
  \[
  \dim \{\, z\in N\mid \dim \xi\inverse(z)\geq r \,\} < \dim N -2r
  \]
  and
  \[
  \dim \{\, z\in N\mid \dim \mu\inverse(z)\geq r \,\} < \dim N -2r.
  \] 
\item[D4] The morphisms $\iM$, $\iP$, and $\iN$ are open embeddings.
\item[D5] The morphisms $\jM$, $\jP$, and $\jN$ are closed embeddings.
\item[D6] A finite group, $\Sigma$, acts on $M_r$ on the right so that
  $N_r\cong M_r/\Sigma$ and $\mu_r$ may be identified with the orbit
  map.
\item[D7] There is a subgroup, $\Sigma'$, of $\Sigma$, so that
  $P_r\cong M_r/\Sigma'$ and $\eta_r$ may be identified with the orbit
  map.
\end{itemize}

Since $\eta$ and $\xi$ are proper morphisms and the squares in the
basic commutative diagram are cartesian, it follows that all the
horizontal maps in the basic commutative diagram are proper morphisms
and that $\mu$, $\mu_r$, and $\mu_0$ are proper morphisms. Thus, if
$f$ is any of the morphisms in the basic commutative diagram except
$\iM$, $\iP$, or $\iN$, then $Rf_*= Rf_!$. Since $\iM$, $\iP$, and
$\iN$ are open embeddings, we have $\iM^*= \iM^!$, $\iP^*= \iP^!$, and
$\iN^*= \iN^!$. Finally, since $\eta_r$, $\xi_r$, and $\mu_r$ are
finite covering maps, we have $\eta_r^!= \eta_r^*$, $\xi_r^!=
\xi_r^*$, $\mu_r^!= \mu_r^*$, $R(\eta_r)_!= (\eta_r)_!$, $R(\xi_r)_!=
(\xi_r)_!$, and $R(\mu_r)_!= (\mu_r)_!$.

In this section we prove the following theorem.

\begin{theorem}\label{thm3.5}
  The group $\Sigma$ acts on $H_\bullet(M_0)$ and there is an
  isomorphism $h'\colon H_{\bullet}(P_0) \cong
  H_{\bullet}(M_0)^{\Sigma'}$ so that if $\Av\colon H_\bullet(M_0) \to
  H_\bullet(M_0) ^{\Sigma'}$ is the averaging map, then the diagram
  \[
     \xymatrix{
       H_{\bullet}(M_0)\ar[rr]^{(\eta_0)_*} \ar[dr]_(.4)\Av & &
       H_{\bullet}(P_0) \ar[dl]^(.4){h'} _(.5){\simeq} \\   
       & H_{\bullet}(M_0)^{\Sigma'} &
     }
   \]
   of graded vector spaces commutes.
\end{theorem}

The idea of the argument is a standard one and is given in the next
three subsections. In \S\ref{3.2} we prove Proposition \ref{prop3.2},
the analog of Theorem \ref{thm3.5} for local systems on $M_r$, $P_r$,
and $N_r$. In \S\ref{3.3} we apply $\IC$ and use that $\xi$ and $\mu$
are small maps to identify the intersection complexes with higher
direct images of constant sheaves. Thus we obtain a sheaf-theoretic
version of Theorem \ref{thm3.5} for complexes of sheaves in $\Dbc(N)$.
In \S\ref{3.4} we complete the proof of the theorem by restricting to
$N_0$, applying $\Ext^\bullet_{N_0} (\BBQ_{N_0}, \jN^!( \,\cdot \,))$,
and showing that the induced map in Borel-Moore homology is
$(\eta_0)_*$.  Since we are concerned not only with complexes of
sheaves, but also the precise maps between them, most of the work
involved is in keeping track of morphisms as we apply the various
functors.

Finally, in \S\ref{3.5} we discuss a two variable version of Theorem
\ref{thm3.5}. Here $M$, $P$, and $N$ are replaced by $M\times M$,
$P\times Q$, and $N\times N$ respectively, $M_0$ and $P_0$ are
replaced by the fibred products $Z =(M\times M) \times_{N\times N}
N_0$ and $X= (P\times Q) \times_{N\times N} N_0$ respectively, and
$\jN$ is replaced by $\delta \jN\colon N_0\to N\times N$, where
$\delta$ is the diagonal map. In the application we are mainly
interested in (see (\ref{eq:bdZ})), $M\times M=\fgt\times \fgt$, $Z$
is the Steinberg variety of $G$, and $X$ is the generalized Steinberg
variety $\XPQ$.

As we have observed above, all the horizontal maps in the basic
commutative diagram are proper, so direct image and direct image with
proper support are the same functors for these maps. Direct image with
proper support is better adapted to Borel-Moore homology, so the
following argument is phrased in terms of direct image with proper
support.

\subsection{}\label{3.2}
First, $\mu_r$ may be identified with the orbit map from $M_r$ to
$M_r/\Sigma$ and so $\Sigma$ acts as automorphisms on the local system
$(\mu_r)_! \BBQ_{M_r}$ on $N_r$. Similarly, $\Sigma'$ acts as
automorphisms on the local system $(\eta_r)_!  \BBQ_{ M_r}$ on $P_r$.

Next, local systems on $N_r$ form an abelian category so we may
consider the $\Sigma'$-invariants of the local system $(\mu_r)_!
\BBQ_{M_r}$.  Let 
\[
\xymatrix{\Av\colon (\mu_r)_! \BBQ_{M_r} \ar[r]& \left((\mu_r)_!
  \BBQ_{M_r}\right) ^{\Sigma'}}
\]
denote the projection onto the local system of $\Sigma'$-invariants
given by averaging over $\Sigma'$.

Finally, recall from \S\ref{2.2} that $\alpha_{\eta_r}\colon \eta_r^*
\BBQ_{P_r} \to \BBQ_{M_r}$ is the natural isomorphism.  Since
$\eta_r^*= \eta_r^!$, we may consider $\alpha_{\eta_r}$ as a map from
$\eta_r^! \BBQ_{P_r}$ to $\BBQ_{M_r}$ and so we may apply
$\Phi_{\eta_r} \inverse$ to $\alpha_{\eta_r} \inverse$ and get a map
from $(\eta_r)_!\BBQ_{M_r}$ to $\BBQ_{P_r}$.  Define
\[
\xymatrix{\gamma_r\colon (\eta_r)_!\BBQ_{M_r} \ar[r]& \BBQ_{P_r}}
\quad \text{by} \quad \gamma_r= \Phi_{\eta_r}\inverse (\alpha_{\eta_r}
\inverse)= \epsilon^{\eta_r} \circ (\eta_r)_!( \alpha_{\eta_r}
\inverse) .
\]

The following proposition is easily proved either directly, or by
using the correspondence between local systems and representations of
fundamental groups.

\begin{proposition}\label{prop3.2}
  There is an isomorphism $h_r\colon (\xi_r)_! \BBQ_{P_r} \cong
  \left((\mu_r)_! \BBQ_{M_r}\right) ^{\Sigma'}$ so that the diagram
  \begin{equation*}
     \xymatrix{
       (\mu_r)_! \BBQ_{M_r} \ar[rr]^{(\xi_r)_! ( \gamma_r)}
       \ar[dr]_(.4)\Av & & (\xi_r)_! \BBQ_{P_r}
       \ar[dl]^(.4){h_r}_(.5){\simeq} \\ 
       & \left((\mu_r)_! \BBQ_{M_r}\right) ^{\Sigma'} &
     }
   \end{equation*}
   of local systems on $N_r$ commutes.
\end{proposition}

\subsection{}\label{3.3}
In this subsection we prove the following proposition, the analog of
Proposition \ref{prop3.2} for $M$, $P$, and $N$.

\begin{proposition}\label{prop3.3}
  There is a map $\gamma \colon R\eta_! \BBQ_{M} \to \BBQ_{P}$ and an
  isomorphism $h\colon R\xi_! \BBQ_{P} \to \left(R\mu_!
    \BBQ_{M}\right) ^{\Sigma'}$ so that the diagram
  \begin{equation*}
     \xymatrix{
       R\mu_! \BBQ_{M} \ar[rr]^{R\xi_! (\gamma )} \ar[dr]_(.4)\Av
       & & R\xi_! \BBQ_{P} \ar[dl]^(.4){h} _(.5){\simeq} \\  
       & \left(R\mu_! \BBQ_{M}\right) ^{\Sigma'} &
     }
   \end{equation*}
   of complexes in $\Dbc(N)$ commutes.
\end{proposition}

We can apply the functor $\IC(N, \,\cdot \,)$ to the diagram of local
systems in Proposition \ref{prop3.2} and obtain a commutative triangle
of complexes in $\Dbc(N)$. Since the functor $\IC(N, \,\cdot \,)$
takes its values in an abelian category of perverse sheaves on $N$ and
is an additive functor by construction, we may consider $\IC(N,
\,\cdot \,)$ as an additive functor between abelian categories. It
follows that $\Sigma$ acts on $\IC(N, (\mu_r)_!  \BBQ_{M_r})$, that
\[
\IC(N, ((\mu_r)_! \BBQ_{M_r})^{\Sigma'}) \cong \IC(N, (\mu_r)_!
\BBQ_{M_r})^{\Sigma'},
\]
and that if $\Av$ is the averaging map, the diagram
\begin{equation*}
   \xymatrix{
   \IC(N, (\mu_r)_! \BBQ_{M_r}) \ar[rr]^{\IC(N, (\xi_r)_! (\gamma_r))}
    \ar[dr]_\Av && \IC(N, (\xi_r)_! \BBQ_{P_r}) \ar[dl]^{\IC(h_r)}
    _(.5){\simeq} \\   
       & \IC(N, (\mu_r)_! \BBQ_{M_r}) ^{\Sigma'} &
     }
\end{equation*}
of complexes in $\Dbc(N)$ commutes.

Since $\xi$ and $\mu$ are small maps, it follows from the axioms
characterizing intersection complexes (see \cite[\S4.13]
{borel:intersection}) that $\IC(N, (\mu_r)_! \BBQ_{M_r})$ and $\IC(N,
(\xi_r)_! \BBQ_{P_r})$ are isomorphic in $\Dbc(N)$ to the direct
images $R\mu_! \BBQ_M$ and $R\xi_! \BBQ_P$ respectively. Moreover,
since the $\Sigma$-action on $R\mu_! \BBQ_M$ comes from transport of
structure from $(\mu_r)_! \BBQ_{M_r}$ it follows that there are
isomorphisms, $\mubar$, $\xibar$, and $h$, so that if $g= \xibar
\inverse \circ \IC(N, (\xi_r)_! (\gamma_r)) \circ \mubar$, then the
diagram
\begin{equation} \label{eq:di}
  \xymatrix{
  \IC(N, (\mu_r)_! \BBQ_{M_r}) \ar[rr]^{\IC(N, (\xi_r)_! (\gamma_r))}
    && \IC(N, (\xi_r)_! \BBQ_{P_r})  \\
   R\mu_! \BBQ_{M}  \ar[u]^\mubar \ar[rr]^{g} \ar[dr]_(.4)\Av &
   & R\xi_! \BBQ_{P} \ar[dl]^(.4){h} \ar[u]_\xibar \\ &
   \left(R\mu_! \BBQ_{M}\right) ^{\Sigma'} & 
   }
\end{equation}
in $\Dbc(N)$ commutes. We can apply the functor $\Ext^\bullet_{N_0}
(\BBQ_{N_0}, \jN^!( \,\cdot \,))$ to the bottom triangle in
(\ref{eq:di}) and obtain a commutative triangle of $\Ext$-groups that
are isomorphic to the Borel-Moore homology groups in the statement of
Theorem \ref{thm3.5}. In order to show that the resulting horizontal
map is indeed the direct image map in Borel-Moore homology induced by
$\eta_0$, we need to choose the isomorphisms $\mubar$ and $\xibar$
appropriately and identify the map $g$ in (\ref{eq:di}). This is
accomplished in the next lemma and the following corollary.

Since $P$ is a purely $d$-dimensional, rational homology manifold, we
have $\BBD_P\cong \BBQ_P[2d]$ in $\Dbc(P)$. We denote by $\nuP$ a
fixed isomorphism, $\nuP\colon \BBD_P\to \BBQ_P[2d]$ in $\Dbc(P)$.

Now $\iM^!\BBD_M[-2d]$ and $\iM^!\BBQ_M$ are in $\Dbc(M_r)$ and
$\iM^!( \alpha_\eta) \circ \iM^!( \eta^!(\nuP) \circ \beta_\eta)$ is
an isomorphism between them, so $\iM^!\BBD_M[-2d]$ is in fact a local
system on $M_r$.  Notice that $\eta^!(\nuP)\colon \eta^! \BBD_P
[-2d]\to \eta^!\BBQ_P$ and $\alpha_\eta\colon \eta^* \BBQ_P\to
\BBQ_M$, so the composition $\alpha_\eta \circ \eta^!(\nuP)$ is not
defined.  However,
\[
\iM^! \eta^!= (\eta\iM)^! =(\iP \eta_r)^!= \eta_r^! \iP^!= \eta_r^*
\iP^*= \iM^* \eta^*= \iM^! \eta^*,
\] 
so the composition $\iM^!(\alpha_\eta) \circ \iM^!( \eta^!(\nuP))$ is
defined.

By \cite[Lemma 4.11]{borel:intersection} there is a unique isomorphism
of local systems on $M$ that restricts to $\iM^!(\alpha_\eta) \circ
\iM^!( \eta^!(\nuP) \circ \beta_\eta)$. The statement in
\cite{borel:intersection} assumes that $M$ is a manifold, but the
argument applies when $M$ is a variety that is a rational homology
manifold. Denote this isomorphism by $\nuM^P$, so $\nuM^P\colon
\BBD_M[-2d]\to \BBQ_M$ and
\begin{equation}
  \label{eq:nuM}
  \iM^!(\nuM^P)= \iM^!(\alpha_\eta) \circ \iM^!( \eta^!(\nuP) \circ
  \beta_\eta).
\end{equation}

Define $\gamma\colon R\eta_!\BBQ_M\to \BBQ_P$ by
\[
\gamma= \nuP \circ \Phi_\eta\inverse (\beta_\eta \circ
(\nuM^P)\inverse)= \nuP \circ \epsilon^\eta \circ R\eta_!(
\beta_\eta\circ (\nuM^P)\inverse) = \Phi_\eta\inverse ( \eta^!(
\nuP)\circ \beta_\eta\circ (\nuM^P) \inverse).
\]

\begin{lemma}\label{lem3.3}
  The diagram
  \[ 
  \xymatrix{ \iN^! R\mu_!\BBQ_M \ar[rr]^{\iN^! R\xi_!(\gamma)}
    \ar[d]_{ (\mu_r)_!(\alpha_{\iM}) \circ \bc_{ \mu_r, \iM}}
    && \iN^! R\xi_! \BBQ_P \ar[d]^{(\xi_r)_!(\alpha_{\iP}) \circ \bc_{
        \xi_r, \iP}}\\ 
    (\mu_r)_! \BBQ_{M_r} \ar[rr]_{(\xi_r)_!( \gamma_r)} && (\xi_r)_!
    \BBQ_{P_r} }
  \]
  of complexes in $\Dbc(N_r)$ commutes.
\end{lemma}

\begin{proof}
  Since $\bc_{\mu_r, \iM}= (\xi_r)_!( \bc_{ \eta_r, \iM}) \circ \bc_{
    \xi_r, \iP}$, we need to show that
  \[
  (\xi_r)_!(\alpha_{\iP}) \circ \bc_{ \xi_r, \iP} \circ \iN^!
  R\xi_!(\gamma) = (\xi_r)_! ( \gamma_r ) \circ
  (\mu_r)_!(\alpha_{\iM}) \circ (\xi_r)_!( \bc_{ \eta_r, \iM}) \circ
  \bc_{ \xi_r, \iP}.
  \]
  Using the naturality of the base change morphism $\bc_{ \xi_r, \iP}$
  we see that it is enough to show that
  \[
  (\xi_r)_!(\alpha_{\iP}) \circ (\xi_r)_!\iP^!( \gamma) = (\xi_r)_!
  (\gamma_r) \circ (\mu_r)_!(\alpha_{\iM}) \circ (\xi_r)_!( \bc_{
    \eta_r, \iM}).
  \]

  Since $\gamma_r= \epsilon^{\eta_r} \circ (\eta_r)_!( \alpha_{\eta_r}
  \inverse)$ it's enough to show that
  \[
  \alpha_{\iP} \circ \iP^!( \gamma) = \epsilon^{\eta_r}\circ
  (\eta_r)_!( \alpha_{\eta_r} \inverse \circ \alpha_{\iM}) \circ \bc_{
    \eta_r, \iM}.
  \]
  
  Equivalently, it's enough to show that
  \[
  \iP^!( \gamma) \circ \bc_{ \eta_r, \iM}\inverse =
  \alpha_{\iP}\inverse \circ \epsilon^{\eta_r} \circ (\eta_r)_!(
  \alpha_{\eta_r} \inverse \circ \alpha_{\iM}).
  \]
  
  Finally, $\eta \iM= \iP\eta_r$ and so $\Phi_{\iM} \Phi_\eta=
  \Phi_{\eta_r} \Phi_{\iP}$ and hence $\Phi_{\iP} \Phi_\eta\inverse=
  \Phi_{\eta_r}\inverse \Phi_{\iM}$. Therefore:
  \begin{alignat*}{2}
    \iP^!( \gamma) \circ \bc_{ \eta_r, \iM}\inverse &= \iP^!\left(
      \Phi_\eta\inverse (\eta^!(\nuP) \circ \beta_\eta \circ (\nuM^P)
      \inverse)\right) \circ \Phi_{\iP} (R\eta_!(\epsilon^{\iM}))& \\
    &=\Phi_{\iP}\left( \Phi_\eta\inverse (\eta^!(\nuP) \circ
      \beta_\eta \circ (\nuM^P) \inverse)\circ R\eta_!
      (\epsilon^{\iM}) \right)& \quad\text{(by \ref{eq:ad2})} \\
    &=\Phi_{\iP} \Phi_\eta\inverse \left(\eta^!(\nuP) \circ \beta_\eta
      \circ (\nuM^P) \inverse \circ \epsilon^{\iM}\right) &
    \quad\text{(by \ref{eq:ad3})}\\
    &=\Phi_{\eta_r}\inverse \Phi_{\iM} \left(\eta^!(\nuP) \circ
      \beta_\eta \circ (\nuM^P) \inverse \circ \epsilon^{\iM}\right)
    &\quad \text{($\Phi_{\iP} \Phi_\eta\inverse= \Phi_{\eta_r}
      \inverse \Phi_{\iM}$)} \\
    &=\Phi_{\eta_r}\inverse \left(\iM^!(\eta^!(\nuP) \circ \beta_\eta
      \circ (\nuM^P)\inverse)\right) & \quad\text{(by \ref{eq:ad2})}\\
    &=\Phi_{\eta_r}\inverse \left(\iM^!(\alpha_\eta \inverse)\right) &
    \quad\text{(by \ref{eq:nuM})} \\
    &=\Phi_{\eta_r}\inverse \left(\iM^*(\alpha_\eta \inverse )\right)
    \\ 
    &=\Phi_{\eta_r}\inverse \left(\eta_r^!(\alpha_{\iP} \inverse)
      \circ \alpha_{\eta_r}\inverse \circ \alpha_{\iM}\right) &
    \quad\text{(by \ref{2.2}.\ref{2.2.2})} \\
    &=\alpha_{\iP} \inverse \circ \Phi_{\eta_r}\inverse \left(
      \alpha_{\eta_r}\inverse \circ \alpha_{\iM}\right) &
    \quad\text{(by \ref{eq:ad3})}\\
    &=\alpha_{\iP} \inverse \circ \epsilon^{\eta_r}\circ (\eta_r)_!
    \left( \alpha_{\eta_r}\inverse \circ \alpha_{\iM}\right) &
    \quad\text{(by \ref{eq:ad1})}
  \end{alignat*}
  This completes the proof of the lemma.
\end{proof}

\begin{corollary}\label{cor3.3}
  There are isomorphisms,
  \[
  \xymatrix{\mubar\colon R\mu_! \BBQ_M \ar[r] & \IC(N, (\mu_r)_!
    \BBQ_{M_r})} \quad \text{and} \quad \xymatrix{\xibar\colon R\xi_!
    \BBQ_P \ar[r] & \IC(N, (\xi_r)_!  \BBQ_{P_r}),}
  \]
  so that the diagram 
  \[
  \xymatrix{ R\mu_!\BBQ_M \ar[rrr]^{R\xi_!( \gamma)} \ar[d] _{\mubar}
    &&& R\xi_! \BBQ_P \ar[d]^{\xibar} \\
    \IC(N, (\mu_r)_! \BBQ_{M_r}) \ar[rrr]_{\IC( (\xi_r)_!(\gamma_r))}
    &&& \IC(N, (\xi_r)_!  \BBQ_{P_r}) } 
  \]
  of complexes in $\Dbc(N)$ commutes.
\end{corollary}

\begin{proof}
  We have already observed that since $\xi$ and $\mu$ are small maps,
  the direct images, $R\xi_!  \BBQ_P$ and $R\mu_! \BBQ_M$ are
  isomorphic in $\Dbc(N)$ to $\IC(N, \xi_!  \BBQ_{P_r})$ and $\IC(N,
  \mu_!  \BBQ_{M_r})$ respectively.  Thus, $R\xi_!  \BBQ_P$ and
  $R\mu_! \BBQ_M$ are in the image of $\IC$. It is shown in
  \cite[Theorem 3.5] {goreskymacpherson:intersectionII} that on the
  image of $\IC$, the composition $\IC(N, \,\cdot\,) \circ \iN^*$ is
  naturally equivalent to the identity so there are isomorphisms,
  \[
  \xymatrix{\ic_\mu\colon R\mu_!\BBQ_M \ar[r]^-{\simeq}& \IC(N, \iN^!
    R\mu_!\BBQ_M)} \quad \text{and} \quad \xymatrix{\ic_\xi\colon
    R\xi_!  \BBQ_P \ar[r]^-{\simeq}& \IC(N, \iN^!
    R(\xi_r)_!\BBQ_{P_r}),}
  \]
  in $D(N)$ with $\iN^*( \ic_\mu)= \id$ and $\iN^*( \ic_\xi)= \id$.
  Since $\IC$ is fully faithful, it follows that the diagram
  \[
  \xymatrix{ R\mu_!\BBQ_M \ar[rrr]^{R\xi_!( \gamma)} \ar[d]_{\ic_\mu}
    &&& R\xi_! \BBQ_P \ar[d]^{\ic_\xi} \\
    \IC(N, \iN^! R\mu_!\BBQ_M) \ar[rrr]_{\IC(\iN^! R\xi_!(\gamma))}
    &&& \IC(N, \iN^! R(\xi_r)_!\BBQ_{P_r})}
  \]
  commutes. 
  
  If we apply $\IC$ to the commutative diagram in the lemma we get a
  commutative diagram:
  \[
  \xymatrix{ \IC(N, \iN^! R\mu_!\BBQ_M) \ar[rrr]^{\IC(\iN^!
      R\xi_!(\gamma))} \ar[d]_{\IC((\mu_r)_!(\alpha_{\iM}) \circ \bc_{
        \mu_r, \iM})} &&& \IC(N, \iN^! R(\xi_r)_!\BBQ_{P_r})
    \ar[d]^{ \IC((\xi_r)_!(\alpha_{\iP}) \circ \bc_{\xi_r, \iP})} \\
    \IC(N, (\mu_r)_!  \BBQ_{M_r}) \ar[rrr]_{\IC((\xi_r)_!( \gamma_r))}
    &&&\IC(N, (\xi_r)_!  \BBQ_{P_r}) }
  \]
  
  Therefore, if we define $\mubar= \IC((\mu_r)_!(\alpha_{\iM}) \circ
  \bc_{ \mu_r, \iM}) \circ \ic_\mu$ and $\xibar= \IC((\xi_r)_!
  (\alpha_{\iP}) \circ \bc_{\xi_r, \iP}) \circ \ic_\xi$, the corollary
  follows.
\end{proof}

Since $\mubar\colon R\mu_! \BBQ_M\to \IC(N, (\mu_r)_! \BBQ_{M_r})$ is
an isomorphism, it follows that $\Sigma$ acts on $R\mu_! \BBQ_M$ by
transport of structure and that $\mubar$ induces an isomorphism
between $\Sigma'$-invariants, say $\mubar'\colon (R\mu_!
\BBQ_M)^{\Sigma'} \to \IC(N, (\mu_r)_!  \BBQ_{M_r})^{\Sigma'}$, which
commutes with the respective averaging maps.

Now consider the diagram:
\begin{equation*}
  \xymatrix@= 17pt{ R\mu_!\BBQ_M \ar[r]^(.4){\mubar} \ar[d]_\Av &
    \IC(N, (\mu_r)_!  \BBQ_{M_r}) \ar[rr]^{\IC((\xi_r)_! (\gamma_r))}
  \ar[dr]_(.4)\Av && \IC(N, (\xi_r)_!  \BBQ_{P_r})
  \ar[dl]^(.35){\IC(h_r)}& R\xi_!\BBQ_P \ar[l]_(.35){\xibar}
  \ar@{.>}[d]^h  \\ 
  (R\mu_!\BBQ_M)^{\Sigma'} \ar[rr]_{\mubar'} && \IC(N, (\mu_r)_!
  \BBQ_{M_r}) ^{\Sigma'} && (R\mu_!\BBQ_M)^{\Sigma'}
  \ar[ll]^(.4){\mubar'} }  
\end{equation*}
If $h$ is defined by $h= (\mubar')\inverse \circ \IC(h_r) \circ
\xibar$, then the diagram commutes. By Corollary \ref{cor3.3}, the
composition across the top row is just $R\xi_!(\gamma)$ and so tracing
around the outside of the diagram we see that $h\circ R\xi_!(\gamma)=
\Av$. This completes the proof of Proposition \ref{prop3.3}.

\subsection{}\label{3.4}
In this subsection, we complete the proof of Theorem \ref{thm3.5}.

\begin{lemma}\label{lem3.4.2}
  There are isomorphisms of graded vector spaces,
  \[
  \xymatrix{J'\colon H_{2d-\bullet}(M_0)\ar[r]^-{\simeq}&
    \Ext^\bullet_{N_0}(\BBQ_{N_0}, \jN^! R\mu_!\BBQ_M)}
  \]
  and
  \[
  \xymatrix{J_1'\colon H_{2d-\bullet}(P_0) \ar[r]^-{\simeq}&
    \Ext^\bullet_{N_0}(\BBQ_{N_0}, \jN^! R\xi_!\BBQ_P)}
  \]
  so that the diagram
  \begin{equation*}
    \xymatrix{
    H_{2d-\bullet}(M_0) \ar[rr]^{(\eta_0)_!} \ar[d]_{J'} &&
    H_{2d-\bullet}(P_0) \ar[d]^{J_1'} \\ 
    \Ext_{N_0}^\bullet( \BBQ_{N_0}, \jN^! R\mu_! \BBQ_M)\ar[rr]^{
     (\jN^!R\xi_!(\gamma))_\sharp} \ar[dr]_(.4)\Av &&
    \Ext_{N_0}^\bullet( \BBQ_{N_0}, \jN^! R\xi_! \BBQ_P)
    \ar[dl]^(.4){(\jN^!(h))_\sharp} \\  
    &\Ext_{N_0}^\bullet( \BBQ_{N_0}, \jN^! R\mu_! \BBQ_M)^{\Sigma'}&}
  \end{equation*}
  commutes.
\end{lemma}

Assuming for a moment that the lemma has been proved, we complete the
proof of Theorem \ref{thm3.5} using the argument at the end of
\S\ref{3.3} as follows.

Since $J'$ is an isomorphism, $\Sigma$ acts on $H_{2d-\bullet}(M_0)$
by transport of structure and $J'$ induces an isomorphism between
$\Sigma'$-invariants, say $\Jbar$, which commutes with the respective
averaging maps.

Now consider the diagram
\[
\xymatrix{ H_{2d-\bullet}(M_0) \ar[r]^(.6){J'} \ar[d]_\Av &E_1
  \ar[rr]^{(\jN^!R\xi_!(\gamma))_\sharp} \ar[dr]_(.4)\Av && E_2
  \ar[dl]^(.35){(\jN^!(h))_\sharp}&
  H_{2d-\bullet}(P_0) \ar[l]_(.6){J_1'} \ar@{.>}[d]^{h'}  \\
  H_{2d-\bullet}(M_0)^{\Sigma'} \ar[rr]_(.6){\Jbar} && E_3 &&
  H_{2d-\bullet}(M_0)^{\Sigma'} \ar[ll]^(.6){\Jbar} }
\]
where 
\begin{gather*}
  E_1= \Ext^\bullet_{N_0} (\BBQ_{N_0}, \jN^! R\mu_!\BBQ_M), \quad E_2=
  \Ext^\bullet_{N_0} (\BBQ_{N_0}, \jN^!  R\xi_!\BBQ_P), \quad
  \text{and}\\ 
  E_3= \Ext^\bullet_{N_0} (\BBQ_{N_0}, \jN^!  R\mu_!\BBQ_M)^{\Sigma'}.
\end{gather*}
If $h'$ is defined by $h'= (\Jbar)\inverse \circ (\jN^!(h))_\sharp
\circ J_1'$, then the diagram commutes.  By Lemma \ref{lem3.4.2}, the
composition across the top row is $(\eta_0)_*$ and so tracing around
the outside of the diagram we see that $h'\circ (\eta_0)_*= \Av$. This
proves Theorem \ref{thm3.5}.

It remains to prove Lemma \ref{lem3.4.2}.

First, we apply the functor $\Ext^\bullet_{N_0} (\BBQ_{N_0}, \jN^!(
\,\cdot \,))$ to the diagram in Proposition \ref{prop3.3} and obtain a
commutative triangle of graded vector spaces. Since the functor
$\Ext_{N_0}^\bullet( \BBQ_{N_0}, \jN^!(\,\cdot \,))$ restricted to the
abelian category of perverse sheaves in which $\IC(N, \,\cdot\,)$
takes its values is an additive functor between abelian categories, it
follows that $\Sigma$ acts on $\Ext_{N_0}^\bullet( \BBQ_{N_0}, \jN^!
R\mu_!  \BBQ_M)$, that
\[
\Ext_{N_0}^\bullet( \BBQ_{N_0}, \jN^! (R\mu_!  \BBQ_M)^{\Sigma'})
\cong \Ext_{N_0}^\bullet( \BBQ_{N_0}, \jN^! R\mu_!
\BBQ_M)^{\Sigma'},
\]
and that if $\Av$ is the averaging map, the diagram of graded vector
spaces
\begin{equation*}
   \xymatrix{
     \Ext_{N_0}^\bullet( \BBQ_{N_0}, \jN^! R\mu_! \BBQ_M)\ar[rr]^{
     (\jN^!R\xi_!(\gamma))_\sharp} \ar[dr]_(.4)\Av &&
   \Ext_{N_0}^\bullet( \BBQ_{N_0}, \jN^! R\xi_! \BBQ_P)
   \ar[dl]^(.4){(\jN^!(h))_\sharp} \\  
     &\Ext_{N_0}^\bullet( \BBQ_{N_0}, \jN^! R\mu_!  \BBQ_M)^{\Sigma'}
     &}
\end{equation*}
commutes.

Next, recall that $\gamma= \nuP \circ \epsilon^\eta \circ R\eta_!(
\beta_\eta\circ (\nuM^P)\inverse)$, so using (\ref{eq:ad1}) we get
\[
\xymatrix{\nuP \circ \Phi_\eta\inverse (\beta_\eta)= \gamma \circ
  R\eta_!(\nuM^P) \colon R\eta_!\BBD_M[-2d] \ar[r]& \BBQ_P.}
\]
Applying $\jN^! R\xi_!$ we get $\jN^!R\xi_! (\nuP) \circ \jN^!  R\xi_!
( \Phi_\eta\inverse (\beta_\eta)) = \jN^!R\xi_! (\gamma) \circ  \jN^!
R\mu_!  (\nuM^P)$. This shows that the diagram
\[
\xymatrix{ \Ext^{\bullet}_{N_0}(\BBQ_{N_0}, \jN^! R\mu_!\BBD_M[-2d])
  \ar[rrr]^{ (\jN^!R\xi_!(\Phi_\eta \inverse (\beta_\eta)))_\sharp}
  \ar[d]_{(\jN^!  R\mu_!  (\nuM^P))_\sharp} &&&
  \Ext^{\bullet}_{N_0}(\BBQ_{N_0},
  \jN^!  R\xi_!\BBD_P[-2d]) \ar[d]^{(\jN^!  R\mu_!  (\nu_P))_\sharp}\\
  \Ext^{\bullet}_{N_0}(\BBQ_{N_0}, \jN^! R\mu_!\BBQ_M) \ar[rrr]_{ (
    \jN^!R\xi_! (\gamma)) _\sharp} &&&
  \Ext^{\bullet}_{N_0}(\BBQ_{N_0}, \jN^!  R\xi_!\BBQ_P) }
  \]
commutes.

Finally, we show that there are isomorphisms
\[
\xymatrix{J\colon H_{-\bullet}(M_0)\ar[r]^-{\simeq}&
  \Ext^\bullet_{N_0}(\BBQ_{N_0}, \jN^! R\mu_!\BBD_M)}
\]
and
\[
\xymatrix{J_1\colon H_{-\bullet}(P_0) \ar[r]^-{\simeq}&
  \Ext^\bullet_{N_0}(\BBQ_{N_0}, \jN^! R\xi_!\BBD_P)}
\]
so that the diagram
\begin{equation}
  \label{eq:J}
  \xymatrix{ \Ext^\bullet_{M_0} (\BBQ_{M_0}, \BBD_{M_0}) \ar[rrr]^{
    (\eta_0)_*} \ar[d]_{J} &&& \Ext^\bullet_{P_0} (\BBQ_{P_0},
  \BBD_{P_0})  \ar[d]^{J_1}\\
  \Ext^{\bullet}_{N_0}(\BBQ_{N_0}, \jN^! R\mu_!\BBD_M)
  \ar[rrr]_{(\jN^!R\xi_!(\Phi_\eta \inverse (\beta_\eta)))_\sharp} &&&
  \Ext^{\bullet}_{N_0}(\BBQ_{N_0}, \jN^!  R\xi_!\BBD_P) }  
\end{equation}
commutes. Once this has been done, set $J'=(\jN^!R\mu_!(
\nuM^P))_\sharp \circ J$ and $J_1'= (\jN^!R\xi_!( \nuP))_\sharp \circ
J_1$. Then $J_1' \circ (\eta_0)_! = ( \jN^!R\xi_! (\gamma)) _\sharp
\circ J'$ and so the diagram in the statement of Lemma \ref{lem3.4.2}
commutes as claimed.

Recall that since $\eta_0$ is a proper map, it induces a map in
Borel-Moore homology. If $\Psi_{\eta_0}$ is the adjunction of the
adjoint pair $(\eta_0^*, (R\eta_0)_*)$, then $(\eta_0)_*$ is the
composition,
\begin{alignat*}{2}
  H_{-\bullet}(M_0) &= \Ext^\bullet_{M_0} (\BBQ_{M_0}, \BBD_{M_0}) &
  \\
  &\cong \Ext^\bullet_{M_0} ( \eta_0^* \BBQ_{P_0}, \BBD_{M_0}) &
  \qquad \text{by $\alpha_{\eta_0}^\sharp$}\\
  &\cong \Ext^\bullet_{P_0} ( \BBQ_{P_0}, R(\eta_0)_*\BBD_{M_0}) &
  \qquad \text{by $\Psi_{\eta_0}$}\\
  &\cong \Ext^\bullet_{P_0} ( \BBQ_{P_0}, R(\eta_0)_!
  \eta_0^!\BBD_{P_0}) & \qquad \text{by $\left(
      R(\eta_0)_!(\beta_{\eta_0})
    \right)_\sharp$}\\
  &\longrightarrow \Ext^\bullet_{P_0} ( \BBQ_{P_0}, \BBD_{P_0}) &
  \qquad \text{by $\left(\epsilon^{\eta_0} \right)_\sharp$}\\
  &= H_{-\bullet}(P_0),
\end{alignat*}
so 
\[
(\eta_0)_*= \left(\epsilon^{\eta_0} \circ R(\eta_0)_!( \beta_{\eta_0})
\right)_\sharp \circ \Psi_{\eta_0} \circ \alpha_{\eta_0}^\sharp =
\Phi_{\eta_0}\inverse (\beta_{\eta_0})_\sharp \circ \Psi_{\eta_0}
\circ \alpha_{\eta_0}^\sharp.
\]
  
Now consider the diagram
\[
\xymatrix{ \Ext^{\bullet}_{N_0}(\BBQ_{N_0}, \jN^! R\mu_!\BBD_M)
  \ar[r]^{(\dagger)} \ar[d]_{(*)} & \Ext_{P_0}^\bullet (\BBQ_{P_0},
  R(\eta_0)_* \BBD_{M_0}) \ar[d]_{(**)} && \Ext^\bullet_{M_0}
  (\BBQ_{M_0}, \BBD_{M_0}) \ar[ll]_(.42){\Psi_{\eta_0} \circ
    \alpha_{\eta_0}^\sharp}
  \ar[dll]^{(\eta_0)_*}\\
  \Ext^{\bullet}_{N_0}(\BBQ_{N_0}, \jN^!  R\xi_!\BBD_P)
  \ar[r]_{(\dagger \dagger)} & \Ext^\bullet_{P_0} (\BBQ_{P_0},
  \BBD_{P_0}) && }
\]
where $(*)= \left(\jN^! R\xi_!(\Phi_\eta \inverse(\beta_\eta))
\right)_\sharp$, $(**)= \Phi_{\eta_0} \inverse(
\beta_{\eta_0})_\sharp$, and $(\dagger)$ and $(\dagger \dagger)$ are
given by the compositions
\begin{alignat*}{2}
  \Ext^{\bullet}_{N_0}(\BBQ_{N_0}, \jN^!  R\mu_!\BBD_M) &\cong
  \Ext^{\bullet}_{N_0} (\BBQ_{N_0}, R(\xi_0)_! \jP^! R\eta_!\BBD_M)
  &\qquad \text{by $\left(\bc_{\xi_0, \jP, } \right)_\sharp$}\\
  &\cong \Ext^{ \bullet}_{N_0}(\BBQ_{N_0}, R(\xi_0)_!R(\eta_0)_!
  \jM^!\BBD_M) &\qquad\text{by $\left( R(\xi_0)_!(\bc_{\eta_0, \jM})
    \right)_\sharp$} \\
  &\cong \Ext^{\bullet}_{N_0}(\BBQ_{N_0}, R(\xi_0)_!R(\eta_0)_!
  \BBD_{M_0})& \qquad\text{by $\left( R(\xi_0)_* R(\eta_0)_*
      (\beta_{\jM} \inverse)\right)_\sharp$}\\
  &\cong\Ext^{\bullet}_{P_0}(\xi_0^*\BBQ_{N_0}, R(\eta_0)_!
  \BBD_{M_0})
  &\qquad\text{by $\Psi_{\xi_0} \inverse$}\\
  &\cong \Ext^{\bullet}_{P_0}(\BBQ_{P_0}, R(\eta_0)_! \BBD_{M_0})
  &\qquad\text{by $\left( \alpha_{\xi_0}\inverse \right)_\sharp $}
\end{alignat*}
and
\begin{alignat*}{2}
  \Ext^{\bullet}_{N_0}(\BBQ_{N_0}, \jN^! R\xi_!\BBD_P)&\cong
  \Ext^{\bullet}_{N_0}(\BBQ_{N_0}, R(\xi_0)_! \jP^! \BBD_P) &\qquad
  \text{by $\left(\bc_{\xi_0, \jP, } \right)_\sharp$} \\
  &\cong \Ext^{\bullet}_{N_0}(\BBQ_{N_0}, R(\xi_0)_! \BBD_{P_0}) &
  \qquad\text{by $\left( R(\xi_0)_*(\beta_{\jP} \inverse)
    \right)_\sharp$}\\
  &\cong \Ext^{\bullet}_{P_0}(\xi_0^*\BBQ_{N_0}, \BBD_{P_0})
  &\qquad\text{by $\Psi_{\xi_0} \inverse$}\\
  &\cong \Ext^{\bullet}_{P_0}(\BBQ_{P_0}, \BBD_{P_0}) &\qquad\text{by
    $\left( \alpha_{\xi_0}\inverse \right)_\sharp $}
\end{alignat*}
respectively, so
\[
(\dagger)= (\alpha_{\xi_0}\inverse)^\sharp \circ \Psi_{\xi_0}\inverse
\circ \left( R(\xi_0)_* R(\eta_0)_* (\beta_{\jM} \inverse) \circ
  R(\xi_0)_*( \bc_{\eta_0, \jM}) \circ \bc_{\xi_0, \jP} \right)_\sharp
\]
and
\[
(\dagger \dagger)= (\alpha_{\xi_0}\inverse)^\sharp \circ
\Psi_{\xi_0}\inverse \circ \left( R(\xi_0)_*( \beta_{\jP} \inverse)
  \circ \bc_{\xi_0, \jP} \right)_\sharp .
\]
  
Assume for a moment that $(**) \circ (\dagger)= (\dagger \dagger)
\circ (*)$ and define
\[
\xymatrix{J =(\dagger)\inverse \circ \Psi_{\eta_0} \circ
  \alpha_{\eta_0}^\sharp \colon \Ext^\bullet_{M_0} (\BBQ_{M_0},
  \BBD_{M_0}) \ar[r]& \Ext^{\bullet}_{N_0}(\BBQ_{N_0}, \jN^!
  R\mu_!\BBD_M)}
\]
and
\[
\xymatrix{J_1= (\dagger\dagger) \inverse \colon \Ext^\bullet_{P_0}
  (\BBQ_{P_0}, \BBD_{P_0}) \ar[r]& \Ext^{\bullet}_{N_0} (\BBQ_{N_0},
  \jN^!  R\xi_!\BBD_P).}
\]
Then $J$ and $J_1$ are isomorphisms and
\[ 
J_1 \circ (\eta_0)_* = (\dagger\dagger)\inverse \circ (\eta_0)_* =
(\jN^!R\xi_!(\Phi_\eta \inverse (\beta_\eta)))_\sharp \circ
[(\dagger)\inverse \circ \Psi_{\eta_0} \circ \alpha_{\eta_0}^\sharp]=
(\jN^!R\xi_!(\Phi_\eta \inverse (\beta_\eta)))_\sharp \circ J
\]
so diagram (\ref{eq:J}) commutes as claimed.
  
It remains to show that $(**) \circ (\dagger)= (\dagger \dagger) \circ
(*)$.  Suppose $h$ is in $\Ext^{\bullet}_{N_0}(\BBQ_{N_0}, \jN^!
R\mu_!  \BBD_M)$. Then
\begin{align*}
  ( (**) \circ (\dagger))(h)&= \Phi_{\eta_0}\inverse (\beta_{\eta_0})
  \circ \Psi_{\xi_0} \inverse \left( R(\mu_0)_!(\beta_{\jM} \inverse)
    \circ R(\xi_o)_!(\bc_{\eta_0, \jM}) \circ \bc_{\xi_0, \jP} \circ h
  \right) \circ
  \alpha_{\xi_0}\inverse\\
  &= \Psi_{\xi_0} \inverse \left( R(\xi_0)_!( \Phi_{\eta_0}\inverse
    (\beta_{\eta_0}) \circ R(\eta_0)_!(\beta_{\jM} \inverse) \circ
    \bc_{\eta_0, \jM}) \circ \bc_{\xi_0, \jP} \circ h \right) \circ
  \alpha_{\xi_0}\inverse.
\end{align*}
  
On the other hand, using the naturality of the base change
$\bc_{\xi_0, \jP}$ we have
\begin{align*}
  ((\dagger \dagger)\circ (*))(h)&= \Psi_{\xi_0}\inverse \left(
    R(\xi_0)_!( \beta_{\jP} \inverse) \circ \bc_{\xi_0, \jP} \circ
    \jN^! R\xi_!( \Phi_{\eta}\inverse ( \beta_\eta)) \circ h \right)
  \circ \alpha_{\xi_0} \inverse \\
  &= \Psi_{\xi_0}\inverse \left( R(\xi_0)_!( \beta_{\jP} \inverse)
    \circ R(\xi_0)_!  \jP^!  ( \Phi_{\eta}\inverse ( \beta_\eta))
    \circ \bc_{\xi_0, \jP} \circ h \right) \circ
  \alpha_{\xi_0}\inverse \\
  &= \Psi_{\xi_0}\inverse \left( R(\xi_0)_!( \beta_{\jP} \inverse
    \circ \jP^!  ( \Phi_{\eta}\inverse ( \beta_\eta))) \circ
    \bc_{\xi_0, \jP} \circ h \right) \circ \alpha_{\xi_0}\inverse
\end{align*}
so it is enough to show that
\[
\Phi_{\eta_0}\inverse (\beta_{\eta_0}) \circ R(\eta_0)_!(\beta_{\jM}
\inverse) \circ \bc_{\eta_0, \jM} = \beta_{\jP} \inverse \circ \jP^!
( \Phi_{\eta}\inverse ( \beta_\eta)) .
\]
  
This last equality is easily proved by a computation similar to the
computation in the proof of Lemma \ref{lem3.3} using the definition of
$\bc_{\eta_0, \jM}$ from \S\ref{2.5}; the identities (\ref{eq:ad1}),
(\ref{eq:ad2}), and (\ref{eq:ad3}); the equality $\Phi_{\jP}
\Phi_\eta\inverse= \Phi_{\eta_0} \inverse \Phi_{\jM}$; and
(\ref{2.2}.\ref{2.2.2}).  We omit the details. This completes the
proof of Lemma \ref{lem3.4.2} and Theorem \ref{thm3.5}.

\subsection{}\label{3.5}
>From now on we denote $\eta$ and $\xi$ by $\eta^P$ and $\xi^P$
respectively.

In this subsection we consider the case when we have two
factorizations of $\mu$, $\mu= \xi^P \circ \eta^P= \xi^Q \circ
\eta^Q$, and the spaces $M$ and $N$ in the basic commutative diagram
(\ref{eq:bd}) are replaced by $M\times M$ and $N\times N$
respectively. So, suppose that $Q$ is a purely $d$-dimensional,
rational homology manifold and that in addition to the assumptions
already made concerning the basic commutative diagram, the diagram
\begin{equation*}
  \label{eq:bd3p}
  \xymatrix{
    M \ar[rr]^{\eta^Q}&& Q\ar[rr]^{\xi^Q}&& N\\
    M_r \ar[rr]^{\eta_r^Q} \ar[u]&& Q_r \ar[rr]^{\xi_r^Q} \ar[u]&& N_r
    \ar[u]} 
\end{equation*}
satisfies conditions D1, D2, D3, D4, and D7 with $P$ replaced by $Q$
and $\Sigma'$ replaced by a possibly different subgroup, $\Sigma''$,
of $\Sigma$.

Let $\delta\colon N\to N\times N$ be the diagonal embedding. Then
$\delta \jN\colon N_0\to N\times N$ is a closed embedding. Define $X$
to be the fibred product $(P\times Q) \times_{N\times N} N_0$ and
define $Z$ to be the fibred product $(M\times M) \times_{N\times N}
N_0$. It follows immediately from the definition that a cartesian
product of two small morphisms is again a small morphism. Therefore,
modifying the notation as indicated, the diagram
\begin{equation} \label{eq:bd3}
  \xymatrix{ 
    Z \ar[rr]^-{\eta_0} \ar[d]_{\jZ}&& X \ar[rr]^-{\xi_0}
    \ar[d]_{\jX}&& N_0 \ar[d]_{\jN}\\
    M\times M \ar[rr]^-{\eta^P\times \eta^Q}&& P\times Q
    \ar[rr]^-{\xi^P \times \xi^Q}&& N\times N\\
    M_r \times M_r \ar[rr] \ar[u]&& P_r \times Q_r  \ar[rr] \ar[u]&&
    N_r\times N_r \ar[u]}
\end{equation}
satisfies conditions D1 -- D7 in \S\ref{3.1}. 

We have the following corollary to Theorem \ref{thm3.5}.

\begin{corollary}\label{cor3.5}
  The group $\Sigma\times \Sigma$ acts on the local system
  $(\mu_r\times \mu_r)_! \BBQ_{M_r\times M_r}$. This action induces an
  action of $\Sigma\times \Sigma$ on $R(\mu\times \mu)_!
  \BBQ_{M\times M}$ and hence an action of $\Sigma\times \Sigma$ on
  $H_\bullet(Z)$ by functoriality and transport of structure via the
  isomorphism
  \[
  \xymatrix{J' \colon H_{\bullet}(Z) \ar[r]& \Ext^{4d-\bullet}_{N_0}
    (\BBQ_{N_0}, \jN^!  \delta^!  R(\mu\times \mu)_!\BBQ_{M\times
      M}).}
  \]
  
  There is an isomorphism $h'\colon H_{\bullet}(X) \to H_{\bullet}
  (Z)^{\Sigma'\times \Sigma''}$ so that if $\Av\colon H_\bullet(Z) \to
  H_\bullet(Z) ^{\Sigma'\times \Sigma''}$ is the averaging map, then
  the diagram
  \begin{equation*}\label{cd3.15}
     \xymatrix{
       H_{\bullet}(Z)\ar[rr]^{\eta_*} \ar[dr]_(.4)\Av
       & & H_{\bullet}(X) \ar[dl]^(.4){h'}_{\simeq} \\  
       & H_{\bullet}(Z)^{\Sigma'\times \Sigma''} &
     }
   \end{equation*}
   of graded vector spaces commutes.  
\end{corollary}


\section{Equivariance}

\subsection{}\label{4.1}
In this section we continue the analysis of diagram (\ref{eq:bd3}) and
consider isomorphisms of graded vector spaces from
\cite[\S8.6]{chrissginzburg:representation}
\[
\xymatrix@1{H_\bullet(Z) \ar[r]^-{J} &\Ext_{N_0}^{-\bullet} \left(
    \BBQ_{N_0}, \jN^! \delta^! R(\mu\times \mu)_! \BBD_{M\times M}
  \right) \ar[r]^K & \Ext_{N_0}^{4n-\bullet}( R(\mu_0)_!\BBQ_{M_0},
  R(\mu_0)_! \BBQ_{M_0}) }
\]
where $\dim N_0=2n$, $J$ is as in \S3.4, and $K$ is defined below.
Notice that
\[
\End_{N_0}(R(\mu_0)_!\BBQ_{M_0})= \Ext_{N_0}^{0}(
R(\mu_0)_!\BBQ_{M_0}, R(\mu_0)_! \BBQ_{M_0})\cong H_{4n}(Z).
\]

Recall that $\dim M= \dim N=d$ and define $l=\codim_N N_0=d-2n$.  From
now on, we assume that $M_0$ and $N_0$ are purely $2n$-dimensional,
rational, homology manifolds. We also assume that the fibred products
$X$ and $Z$ in \S\ref{3.5} are purely $2n$-dimensional varieties.

The graded vector space $\Ext_{N_0}^{4n-\bullet}( R(\mu_0)_!
\BBQ_{M_0}, R(\mu_0)_!  \BBQ_{M_0})$ is a graded $\BBQ$-algebra and
the composition $K\circ J$ can be used to give $H_\bullet(Z)$ a
$\BBQ$-algebra structure with $H_i(Z) \cdot H_j(Z) \subseteq
H_{i+j-4n}(Z)$. Since the multiplication in $\Ext_{N_0}^{4n-\bullet}(
R(\mu_0)_!  \BBQ_{M_0}, R(\mu_0)_!  \BBQ_{M_0})$ is composition, we
have
\[
\Ext_{N_0}^{4n-\bullet}( R(\mu_0)_!  \BBQ_{M_0}, R(\mu_0)_!
\BBQ_{M_0}) \cong H_\bullet(Z)\op.
\]

We saw in \S3.3 that $\Sigma$ acts on $R\mu_!\BBQ_{M}$. This action
induces a degree-preserving action of $\Sigma \times \Sigma$ on
$\Ext_{N_0}^{4n-\bullet}( R(\mu_0)_!\BBQ_{M_0}, R(\mu_0)_!
\BBQ_{M_0})$. On the other hand, as in \S3.5, $\Sigma \times \Sigma$
acts on $R(\mu\times \mu)_!  \BBQ_{M\times M}$. This action induces a
degree-preserving $\Sigma\times \Sigma$-action on
$\Ext_{N_0}^{-\bullet} \left( \BBQ_{N_0}, \jN^!  \delta^!  R(\mu\times
  \mu)_! \BBD_{M\times M} \right)$.

In this section we show that the isomorphisms $J$ and $K$ are
$\Sigma\times \Sigma$-equivariant. It then follows that if
$\Sigma\times \Sigma$ acts on the group algebra $\BBQ\Sigma$ in the
usual way, then there are $\Sigma \times \Sigma$-equivariant,
$\BBQ$-algebra homomorphisms
\begin{equation}
  \label{eq:hom}
  \xymatrix{\BBQ \Sigma \ar[r]& \End_{N_0}( R(\mu_0)_! \BBQ_{M_0})
  \ar[r]^-{\simeq}& H_{4n}(Z)\op.} 
\end{equation}

In \S4.2 we describe the $\Sigma \times \Sigma$-action on
$\Ext_{N_0}^{4n-\bullet}( R(\mu_0)_!\BBQ_{M_0}, R(\mu_0)_!
\BBQ_{M_0})$. In \S4.3 we describe the $\Sigma \times \Sigma$-action
on $\Ext_{N_0}^{ -\bullet} \left( \BBQ_{N_0}, \jN^!  \delta^!
  R(\mu\times \mu)_! \BBD_{M\times M} \right)$ and observe that $J$ is
$\Sigma \times \Sigma$-equivariant.  In \S4.4 we define the map $K$,
and in \S4.5--\S4.8 we show that $K$ is $\Sigma \times
\Sigma$-equivariant.

\subsection{}\label{4.2}
We first consider the $\Sigma \times \Sigma$-action on $\Ext_{N_0}^{
  \bullet}( R(\mu_0)_!\BBQ_{M_0}, R(\mu_0)_!  \BBQ_{M_0})$.  Returning
to our original basic commutative diagram (\ref{eq:bd}), $\Sigma$ acts
on the direct image, $(\mu_r)_!  \BBQ_{M_r}$. This action induces a
$\BBQ$-algebra homomorphism
\[
\xymatrix{L_r\colon \BBQ \Sigma\ar[r]& \End_{N_r}
  \left((\mu_r)_!\BBQ_{M_r} \right).}
\]
Applying $\IC$ and transporting the action via the isomorphism $\mubar
\colon R\mu_!  \BBQ_M\to \IC(N, (\mu_r)_! \BBQ_{M_r})$ from Corollary
\ref{cor3.3} gives rise to a $\BBQ$-algebra homomorphism
\[
\xymatrix{L\colon \BBQ \Sigma\ar[r]& \End_N (R\mu_!\BBQ_M)}
\]
with $L(\sigma)= \mubar\inverse \circ \IC(L_r(\sigma)) \circ \mubar$.

Since $L$ is a ring homomorphism, we get an action of $\Sigma \times
\Sigma$ on  $\End_N (R\mu_!\BBQ_M)$ with
\[
(\sigma, \sigma')\cdot f= L(\sigma') \circ f\circ L(\sigma \inverse)
\]
for $f$ in $ \End_N (R\mu_!\BBQ_M)$. 

Clearly, if $\Sigma\times \Sigma$ acts on $\BBQ \Sigma$ by $(\sigma,
\sigma') \cdot x= \sigma' x \sigma\inverse$, then $L$ is $\Sigma
\times \Sigma$-equivariant.

Let $\bc^*\colon \jN^* R\mu_!\to R(\mu_0)_! \jM^*$ be as in \S2.5.
Then $R(\mu_0)_!(\alpha_{\jM}) \circ \bc^*$ is an isomorphism between
$\jN^* R\mu_!  \BBQ_M$ and $R(\mu_0)_!  \BBQ_{M_0}$. We define
\[
\xymatrix{L_0\colon \BBQ\Sigma \ar[r]& \End_{N_0} \left( R(\mu_0)_!
    \BBQ_{M_0}\right)}
\] 
by
\[
L_0(\sigma)= R(\mu_0)_! (\alpha_{\jM}) \circ \bc^* \circ
\jN^*L(\sigma) \circ (\bc^*)\inverse \circ R(\mu_0)_!
(\alpha_{\jM}\inverse).
\]

Since $\Ext^{j}_{N_0}( R(\mu_0)_!\BBQ_{M_0}, R(\mu_0)_! \BBQ_{M_0})
=\Hom_{N_0}( R(\mu_0)_! \BBQ_{M_0}, R(\mu_0)_! \BBQ_{M_0}[j])$ is
naturally an $\End_{N_0} (R(\mu_0)_! \BBQ_{M_0})$-bimodule, we may
define an action of $\Sigma \times \Sigma$ on the graded vector space
$\Ext^{\bullet}_{N_0}( R(\mu_0)_!\BBQ_{M_0}, R(\mu_0)_! \BBQ_{M_0})$
by
\[
(\sigma, \sigma')\cdot g= L_0(\sigma') \circ g\circ L_0(\sigma
\inverse)
\]
for $\sigma$ and $\sigma'$ in $\Sigma$ and $g$ in $\Ext^{
  \bullet}_{N_0}( R(\mu_0)_!\BBQ_{M_0}, R(\mu_0)_! \BBQ_{M_0})$.

\subsection{}\label{4.3}
Next we consider the $\Sigma \times \Sigma$-action on
$\Ext_{N_0}^{\bullet} \left( \BBQ_{N_0}, \jN^!  \delta^!  R(\mu\times
  \mu)_! \BBD_{M\times M} \right)$. Since $M$ is a rational homology
manifold, so is $M\times M$ and we denote by $\nu_{M\times M}$ a fixed
isomorphism, $\nu_{M\times M} \colon \BBD_{M\times M}\to \BBQ_{M\times
  M}[4d]$.

As in \S\ref{3.4} and \S\ref{3.5}, $\Sigma\times \Sigma$ acts as
automorphisms on $R(\mu\times \mu)_!\BBQ_{M\times M}$ and we transport
the group action on $R(\mu\times \mu)_!\BBQ_{M\times M}$ to an action
on $R(\mu\times \mu)_!\BBD_{M\times M}$ using $R(\mu\times \mu)_!
(\nu_{M\times M})$. The group actions induce ring homomorphisms
\[ 
\xymatrix{L_2\colon \BBQ(\Sigma\times \Sigma)\ar[r]& \End_{N\times
    N}(R(\mu\times \mu)_!  \BBQ_{M\times M})}
\]
and
\[
\xymatrix{L_2'\colon \BBQ(\Sigma \times \Sigma)\ar[r]& \End_{N\times
    N}(R( \mu\times \mu)_!\BBD_{M \times M})}
\] 
where $L_2$ and $L_2'$ are related by
\[
L_2'(\sigma, \sigma')= R(\mu\times \mu)_! (\nu_{M\times M} \inverse)
\circ L_2(\sigma, \sigma') \circ R(\mu\times \mu)_!  (\nu_{M\times
  M}).
\]

Notice that $L_2'$ depends on the choice of the orientation
$\nu_{M\times M}$. 

Applying $\Ext_{N_0}^\bullet (\BBQ_{N_0}, \jN^!\delta^!(\, \cdot\,))$
to $R(\mu\times \mu)_!  \BBD_{M\times M}$ and using $L_2'$ we get an
action of $\Sigma\times \Sigma$ on $\Ext_{N_0 }^{\bullet} (\BBQ_{N_0},
\jN^!  \delta^! R(\mu\times \mu)_!  \BBD_{M\times M})$ with
\[
(\sigma, \sigma')\cdot f= (\jN^!\delta^!L_2'( \sigma,
\sigma'))_\sharp (f)= \jN^!\delta^!L_2'( \sigma, \sigma') \circ f
\]
for $f$ in $\Ext_{N_0}^{\bullet} (\BBQ_{N_0}, \jN^!  \delta^!
R(\mu\times \mu)_!  \BBD_{M\times M})$.

As in \S\ref{3.4} and \S\ref{3.5}, the $\Sigma\times \Sigma$-action on
$\Ext_{N_0}^{\bullet} (\BBQ_{N_0}, \jN^!  \delta^! R(\mu\times \mu)_!
\BBQ_{M\times M})$ induces an action of $\Sigma \times \Sigma$ on
$H_\bullet(Z)$ by transport of structure using the isomorphism
\[
\xymatrix{J'=(\jN^! \delta^! R(\mu \times \mu)_!( \nu_{M\times
    M}))_\sharp \circ J \colon H_{4d-\bullet}(Z) \ar[r]&
  \Ext^\bullet_{N_0} (\BBQ_{N_0}, \jN^!  \delta^!  R(\mu\times
  \mu)_!\BBQ_{M\times M}).}
\]

It follows from the definitions that $(\jN^! \delta^! R(\mu \times
\mu)_!( \nu_{M\times M}))_\sharp$ is $\Sigma \times
\Sigma$-equivariant.  This proves the following proposition.

\begin{proposition}\label{prop4.3.1}
  The isomorphism 
  \[
  \xymatrix{J\colon H_\bullet(Z)\ar[r]& \Ext_{N_0}^{-\bullet}
    (\BBQ_{N_0}, \jN^!  \delta^! R(\mu\times \mu)_! \BBD_{M\times M})}
  \]
  is $\Sigma \times \Sigma$-equivariant.
\end{proposition}

\subsection{}\label{4.4}
Define
\[
\xymatrix{K\colon \Ext_{N_0}^{\bullet} (\BBQ_{N_0}, \jN^!\delta^!
  R(\mu\times \mu)_!  \BBD_{M\times M}) \ar[r]&
  \Ext_{N_0}^{\bullet+4n} ( R(\mu_0)_! \BBQ_{M_0}, R(\mu_0)!
  \BBQ_{M_0})}
\]
to be the composition
\begin{align*}
  \Ext_{N_0}^\bullet( \BBQ_{N_0}, \jN^!\delta^! R(\mu\times \mu)_!
  \BBD_{M\times M}) &\xrightarrow{\jN^! \delta^!(k')_\sharp}
  \Ext_{N_0}^\bullet \left(\BBQ_{N_0}, \jN^! \delta^! (R\mu_!
    \BBD_{M}\boxtimes R\mu_! \BBD_{M})\right) \\
  &\xrightarrow{\jN^! \delta^!(c\inverse \boxtimes \id )_\sharp}
  \Ext_{N_0}^\bullet \left(\BBQ_{N_0}, \jN^!\delta^! ((R\mu_!
    \BBQ_{M})\chk \boxtimes R\mu_! \BBD_{M})\right) \\
  &\xrightarrow{\jN^! \delta^!(\lambda)_\sharp} \Ext_{N_0}^\bullet
  \left(\BBQ_{N_0}, \jN^! \delta^! (R\CHom(p^*R\mu_! \BBQ_{M}, q^!
    R\mu_!  \BBD_{M})\right)  \\
  &\xrightarrow{(\nat_{\delta \jN} )_\sharp} \Ext_{N_0}^\bullet
  \left(\BBQ_{N_0}, R\CHom( \jN^* R\mu_! \BBQ_{M}, \jN^!
    R\mu_!  \BBD_{M})\right)  \\
  &\xrightarrow{\can\inverse} \Ext_{N_0}^\bullet \left( \jN^* R\mu_!
    \BBQ_{M}, \jN^!   R\mu_! \BBD_{M}\right) \\
  &\xrightarrow{(a\inverse)^\sharp \circ b_\sharp} \Ext_{N_0}^{
    \bullet+4n} \left( R(\mu_0)_!  \BBQ_{M_0}, R(\mu_0)_!
    \BBQ_{M_0}\right)
\end{align*}
where the notation is as follows:
\begin{itemize}
\item $k'\colon R(\mu\times \mu)_!  \BBD_{M\times M} \to R\mu_!\BBD_M
  \boxtimes R\mu_!\BBD_M$ is the K\"unneth isomorphism (recall that
  $\mu$ is proper).
\item $c= R\mu_!(\dc_M\inverse \circ (\beta_\mu\inverse)_\sharp)
  )\circ \phi_\mu \colon \left( R\mu_! \BBQ_M \right) \chk \to R\mu_!
  \BBD_M$ where $\dc_M$ is as in \S\ref{2.1} and $\beta_\mu$ and
  $\phi_\mu$ are as in \S\ref{2.2}.  Notice that $c$ is an isomorphism
  in $\Dbc(N)$, so $\jN^! \delta^! (c\inverse \boxtimes \id)_\sharp$
  makes sense.
\item $\lambda$, $\nat_\delta$, and $\can$ are as in \S2.
\item $a= R(\mu_0)_!( \alpha_{\jM}) \circ \bc^*\colon \jN^* R\mu_!
  \BBQ_M \to R(\mu_0)_!\BBQ_{M_0}$ (see \S\ref{4.1}).
\item $b= R(\mu_0)_!  (\nu_{M_0} \circ \beta_{\jM}\inverse) \circ
  \bc^! \colon \jN^! R\mu_! \BBD_M \to R(\mu_0)_!\BBQ_{M_0}$ where
  $\bc^!\colon \jN^!R\mu_!\to R(\mu_0)_! \jM^!$ is as in \S\ref{2.5},
  $\beta_{\jM}$ is as in \S\ref{2.2}, and $\nu_{M_0}\colon \BBD_{M_0}
  \to \BBQ_{M_0}[4n]$ is an isomorphism in $\Dbc(M_0)$ (recall that
  $M_0$ is a rational homology manifold).
\end{itemize}

Since $K$ is a composition of isomorphisms of graded vector spaces, it
follows that $K$ is an isomorphism of graded vector spaces that
increases the grading by $4n$. 

\begin{theorem}\label{thm4.4}
  The isomorphism
  \[
  \xymatrix{K\colon \Ext_{N_0}^{\bullet} (\BBQ_{N_0}, \jN^! \delta^!
    R(\mu\times \mu)_!  \BBD_{M\times M}) \ar[r]&
    \Ext_{N_0}^{\bullet+4n} ( R(\mu_0)_!  \BBQ_{M_0}, R(\mu_0)_!
    \BBQ_{M_0})}
  \]
  is $\Sigma\times \Sigma$-equivariant.
\end{theorem}

To prove the theorem we show that $\jN^! \delta^!(k')_\sharp$, $\jN^!
\delta^!(c\inverse \boxtimes \id )_\sharp$, $\can\inverse \circ
(\nat_{\delta \jN} \circ \jN^!  \delta^!  (\lambda))_\sharp$, and
$(a\inverse)^\sharp \circ b_\sharp$ are $\Sigma \times
\Sigma$-equivariant in \S4.5, \S4.6, \S4.7, and \S4.8 respectively.

\subsection{}\label{4.5}
In the situation of \S\ref{3.5} we have two factorizations of $\mu$:
$\mu= \xi^P \eta^P= \xi^Q \eta^Q$. Let $\nuM^P$ and $\nuM^Q$ be two
isomorphisms, $\BBD_M \xrightarrow{\simeq} \BBQ_M[2d]$.  Then $\nuM^P
\boxtimes \nuM^Q \colon \BBD_M \boxtimes \BBD_M \to \BBQ_M \boxtimes
\BBQ_M[4d]$ is an isomorphism in $\Dbc(M\times M)$.  The superscripts
$P$ and $Q$ do not necessarily have anything to do with $P$ and $Q$,
but are convenient for distinguishing between the factors.

Using the orientations $\nuM^P$ and $\nuM^Q$ we can define
$\BBQ$-algebra homomorphisms
\[
\xymatrix{L'_P\colon \BBQ \Sigma \ar[r]& \End_N (R\mu_!\BBD_M)}\qquad
\text{and} \qquad \xymatrix{L'_Q\colon \BBQ \Sigma \ar[r]& \End_N
  (R\mu_!\BBD_M)}
\]
by $L'_P(\sigma)= R\mu_!(\nuM^P) \inverse\circ L(\sigma) \circ
R\mu_!(\nuM^P)$ and $L'_Q(\sigma)= R\mu_!(\nuM^Q) \inverse\circ
L(\sigma) \circ R\mu_!(\nuM^Q)$ respectively.

In the following, we always assume that $\nu_{M\times M}$ is chosen
so that 
\[
\nu_{M\times M}= (k'')\inverse \circ (\nuM^P \boxtimes \nuM^Q) \circ
k'
\]
where 
\[
\xymatrix{k'\colon R(\mu\times \mu)_!  \BBD_{M\times M}
  \ar[r]^-{\simeq}& R\mu_!\BBD_M \boxtimes R\mu_!\BBD_M}
\]
and 
\[
\xymatrix{k''\colon R(\mu\times \mu)_!  \BBQ_{M\times M}
  \ar[r]^-{\simeq}& R\mu_!\BBQ_M \boxtimes R\mu_!\BBQ_M}
\]
are K\"unneth isomorphisms.

The next lemma follows from the naturality of $k'$.

\begin{lemma}\label{lem4.3}
  For $\sigma$ and $\sigma'$ in $\Sigma$, the diagram
  \[
  \xymatrix{ R(\mu\times \mu)_! \BBD_{M\times M} \ar[rr]^{k'} \ar[d]_{
      L_2'(\sigma, \sigma')} && R\mu_! \BBD_M \boxtimes R\mu_! \BBD_M
    \ar[d]^{L_P'(\sigma) \boxtimes L_Q'(\sigma')} \\
    R(\mu\times \mu)_! \BBD_{M\times M} \ar[rr]^{k'} && R\mu_! \BBD_M
    \boxtimes R\mu_! \BBD_M }
  \]
  commutes.
\end{lemma}

The lemma shows that if $\Sigma \times \Sigma$ acts on
$\Ext_{N_0}^\bullet \left(\BBQ_{N_0}, \jN^! \delta^! (R\mu_!
  \BBD_{M}\boxtimes R\mu_! \BBD_{M})\right)$ by
\[
(\sigma, \sigma')\cdot f= (\jN^!\delta^! (L_P'(\sigma) \boxtimes
L_Q'(\sigma')) \circ f,
\]
then $\jN^! \delta^!(k')_\sharp$ is $\Sigma \times \Sigma$-equivariant.

\subsection{}\label{4.6}
In this subsection we show that if $\Sigma \times \Sigma$ acts on
$\Ext_{N_0}^\bullet \left(\BBQ_{N_0}, \jN^!\delta^!  ((R\mu_!
  \BBQ_{M})\chk \boxtimes R\mu_! \BBD_{M})\right)$ by
\[
(\sigma, \sigma')\cdot f= (\jN^!\delta^! (L(\sigma\inverse)\chk
\boxtimes L_Q'(\sigma')) \circ f,
\]
then $\jN^! \delta^!(c\inverse \boxtimes \id )_\sharp$ is
$\Sigma\times \Sigma$-equivariant. In order to do this, it is enough
to show that $c\colon \left( R\mu_! \BBQ_M \right) \chk \to R\mu_!
\BBD_M$ intertwines $L(\sigma\inverse)\chk$ and $L_P'(\sigma)$ for
$\sigma$ in $\Sigma$.

In the rest of this subsection, we denote $\nuM^P$ and $L_P'$ simply
by $\nuM$ and $L'$ respectively.

It is shown in \cite[Theorem 9.8]{borel:intersection} that the Verdier
dual of the intersection complex of a local system is, up to a shift,
the intersection complex of the dual local system. Also, in the
equivalence between local systems and representations of the
fundamental group, the dual of a local system corresponds to the
contragredient representation and the direct image of local systems
corresponds to the induced representation. On the representation
theory side, we are considering permutation representations, which are
obviously equivalent to their contragredients, so it is natural to
expect that for $\sigma$ in $\Sigma$, the Verdier dual of $\sigma$,
acting on $(R\mu_! \BBQ_M) \chk$, may be identified with
$\sigma\inverse$ acting on $(R\mu_! \BBQ_M)$. This is indeed the case
and the next proposition gives the precise formulation we need.

\begin{proposition}\label{prop4.2}
  If $c= R\mu_!(\dc_M\inverse \circ (\beta_\mu\inverse)_\sharp) )\circ
  \phi_\mu$, then the diagram
  \[
  \xymatrix{ \left( R\mu_! \BBQ_M \right)\chk
    \ar[d]_{c} \ar[rr]^{
      L(\sigma\inverse) \chk} && \left( R\mu_! \BBQ_M \right)\chk
    \ar[d]^{c} \\ 
    R\mu_!\BBD_M \ar[rr]_{L'(\sigma)} && R\mu_!\BBD_M }
  \]
  of isomorphisms of complexes in $\Dbc(N)$ commutes for every
  $\sigma$ in $\Sigma$.
\end{proposition}

\begin{proof}
  It follows from \cite[Theorem 9.8]{borel:intersection} that there is
  a unique isomorphism,
  \[
  \xymatrix{\vd\colon \IC(N, (\mu_r)_!\BBQ_{M_r})\chk [-2d] \ar[r]&
    \IC(N, ((\mu_r)_!\BBQ_{M_r})\chk [-2d])}
  \]
  with the property that $\iN^*(\vd)= (\beta_{\iN} \inverse)_\sharp
  \circ \nat_{\iN}$.
  
  Define $\nuMr= \alpha_{\iM} \circ \iM^!(\nuM) \circ \beta_{\iM}$, so
  $\nuMr\colon \BBD_{M_r}\to \BBQ_{M_r}[2d]$ is an isomorphism.
  
  Now consider the ``cube''
  \begin{equation}
  \label{eq:cube}
  \xymatrix{ \IC(N, (\mu_r)_!\BBQ_{M_r})\chk [-2d] \ar[rr]^{
      \IC(L_r(\sigma\inverse )\chk)} \ar[ddd]^{x} \ar[ddddr]^{
      \mubar\chk}  &&\IC(N, (\mu_r)_!\BBQ_{M_r}) \chk  [-2d]
    \ar[ddd]^{x} \ar[ddddr]^{\mubar\chk} &\\
  &&& \\
  &&& \\
  \IC(N, (\mu_r)_!\BBQ_{M_r}) \ar[rr]^{ \IC(L_r( \sigma))}
  \ar[ddddr]^{ \mubar\inverse}  &&\IC(N, (\mu_r)_!\BBQ_{M_r})
    \ar[ddddr]^{\mubar\inverse} &\\   
  & (R\mu_!\BBQ_M)\chk [-2d] \ar[rr]^{L(\sigma\inverse)\chk}
    \ar[ddd]^{y} && (R\mu_!\BBQ_M)\chk [-2d] \ar[ddd]^{y} \\ 
  &&& \\
  &&& \\
  & R\mu_!\BBQ_M \ar[rr]^{L(\sigma)}  && R\mu_!\BBQ_M }
  \end{equation}
  where $x= \IC\left( (\mu_r)_!( \nuMr \circ \dc_{M_r} \inverse \circ
    (\beta_{\mu_r} \inverse)_\sharp ) \circ \phi_{\mu_r} \right) \circ
  \vd$, $y= R\mu_!(\nuM) \circ c$, and $\phi_\mu$ is as in
  \S\ref{2.2}. It follows from the definitions of $y$ and $L'$ that it
  is enough to show that the front face commutes. We show that all
  faces besides the front face commute and so the front face must
  commute also.
  
  The top and bottom faces of (\ref{eq:cube}) commute by definition
  and the left and right faces are equal, so we need to show that the
  back face and the left face commute.
  
  To show that the left face of (\ref{eq:cube}) commutes we need to
  show that the diagram
  \begin{equation}
    \label{eq:dic}
    \xymatrix{ \IC\left(N, (\mu_r)_! \BBQ_{M_r} \right) \chk [-2d]
    \ar[r]^-{\mubar\chk} \ar[d]_{\vd} & \left(R\mu_! \BBQ_{M} \right)
    \chk [-2d] \ar[dd]^{ R\mu_! ((\beta_\mu \inverse)_\sharp) \circ
      \phi_\mu} \\
    \IC\left(N, ((\mu_r)_! \BBQ_{M_r}) \chk [-2d] \right) \ar[d] _{
      \IC ( (\mu_r)_! ((\beta_{\mu_r} \inverse)_\sharp) \circ
      \phi_{\mu_r})}  &\\
    \IC\left(N, (\mu_r)_! \BBQ_{M_r} \chk [-2d] \right)
    \ar@{-->}[r]^-{z} \ar[d]_{\IC ( (\mu_r)_! (\dc_{M_r} \inverse))} &
    R\mu_! \BBQ_{M} \chk [-2d] \ar[d]^{R\mu_!(\dc_M\inverse)} \\
    \IC\left(N, (\mu_r)_! \BBD_{M_r} [-2d] \right) \ar@{-->}[r]^-{w}
    \ar[d]_{\IC ( (\mu_r)_! (\nuMr))} & R\mu_! \BBD_{M} [-2d]
    \ar[d]^{R\mu_!(\nuM)} \\
    \IC\left(N, (\mu_r)_! \BBQ_{M_r} \right) \ar[r]^-{\mubar\inverse}
    &R\mu_! \BBQ_{M} }
  \end{equation}
  commutes.
 
  For the rest of this proof, set $\bc=\bc_{ \mu_r, \iM}$.
 
  As in the proof of Corollary \ref{cor3.3}, since we have
  $R\mu_!\BBQ_M \chk [-2d] \cong \IC\left(N, (\mu_r)_! \BBQ_{M_r} \chk
    [-2d] \right)$ and $R\mu_!\BBD_M [-2d] \cong \IC\left(N, (\mu_r)_!
    \BBD_{M_r} [-2d] \right)$, there are isomorphisms,
  \[
  \xymatrix{\ic_c \colon R\mu_! \BBQ_M \chk [-2d] \ar[r]& \IC\left(N,
      (\mu_r)_!  \BBQ_{M_r} \chk [-2d] \right)}
  \]
  and
  \[
  \xymatrix{\ic_d \colon R\mu_!\BBD_M [-2d] \ar[r]& \IC\left(N,
      (\mu_r)_!  \BBD_{M_r} [-2d] \right),}
  \]
  in $D(N)$ with $\iN^*(\ic_c)= \id$ and $\iN^*( \ic_d)= \id$.  Define
  \[
  z= \ic_c\inverse \circ \IC( \bc\inverse \circ (\mu_r)_!  (\nat_{\iM}
  \inverse \circ (\beta_{\iM})_\sharp \circ \alpha_{\iM}^\sharp )),
  \quad w= \ic_d\inverse \circ \IC( \bc\inverse \circ (\mu_r)_!
  (\beta_{\iM}))
  \] 
  and recall that $\mubar\inverse= \ic_\mu\inverse \circ \IC(
  \bc\inverse \circ (\mu_r)_!  (\alpha_{\iM}\inverse ))$.
  
  Since all the complexes in (\ref{eq:dic}) are in the image of $\IC$,
  it is enough to show that (\ref{eq:dic}) commutes after applying
  $\iN^*$.
  
  First, it follows from the definition of $\nuMr$ and the naturality
  of $\bc$ that
  \begin{align*}
    \iN^* \left( \mubar\inverse \circ \IC((\mu_r)_!(\nuMr)) \right)&=
    \bc\inverse \circ (\mu_r)_!  (\alpha_{\iM}\inverse) \circ
    (\mu_r)_!(\nuMr)\\
    &= \iN^* R\mu_!(\nuM) \circ \bc\inverse \circ (\mu_r)_!
    (\beta_{\iM}) \\
    &= \iN^* \left( R\mu_!(\nuM) \circ w\right).
  \end{align*}
  
  Second, it follows from (\ref{2.2}.\ref{2.2.1}) applied to $\iM$ and
  the naturality of $\bc$ that
  \begin{align*}
    \iN^* \left( R\mu_!(\dc_M\inverse) \circ z\right) &= \iN^*
    R\mu_!(\dc_M\inverse) \circ \bc\inverse \circ (\mu_r)_!
    (\nat_{\iM}
    \inverse \circ (\beta_{\iM})_\sharp \circ \alpha_{\iM}^\sharp) \\
    &= \bc\inverse \circ (\mu_r)_!  (\beta_{\iM} \circ \dc_{M_r}\inverse) \\
    &=\iN^* \left(w \circ \IC((\mu_r)_!( \dc_{M_r}\inverse)) \right).
  \end{align*}
  
  Lastly, it follows from the naturality of $\nat_{\iN}$,
  $\nat_{\iM}$, $\bc$, $\phi_\mu$, and $\phi_{\mu_r}$, Lemma
  \ref{2.5}, (\ref{2.2}.\ref{2.2.2}), and the equality $\iN^* (\vd)=
  (\beta_{\iN} \inverse)_\sharp \circ \nat_{\iN}$ that
  \begin{align*}
    \iN^* \left( R\mu_! ((\beta_\mu \inverse)_\sharp) \circ \phi_\mu
      \circ \mubar\chk \right) &= \iN^* \left( R\mu_! ((\beta_\mu
      \inverse)_\sharp) \circ \phi_\mu \right) \circ \nat_{\iN} \circ
    \inverse \left( (\mu_r)_!  (\alpha_{\iM}) \circ \bc \right)^\sharp
    \circ  \nat_{\iN} \\
    &= \bc\inverse \circ (\mu_r)_!  (\nat_{\iM} \inverse \circ
    (\beta_{\iM})_\sharp \circ \alpha_{\iM}^\sharp \circ
    (\beta_{\mu_r} \inverse)_\sharp) \circ \phi_{\mu_r} \circ
    (\beta_{\iN} \inverse)_\sharp \circ \nat_{\iN} \\
    &= \iN^* \left(z \circ \IC ( (\mu_r)_!  ((\beta_{\mu_r}
      \inverse)_\sharp) \circ \phi_{\mu_r} )\circ \vd \right).
  \end{align*}
  
  Finally, consider the back face of diagram (\ref{eq:cube}). It
  follows from the uniqueness of $\vd$ that $\vd\circ \IC(L(\sigma
  \inverse)) \chk= \IC( L_r( \sigma\inverse) \chk) \circ \vd$. Thus,
  to show that the back face commutes, it is enough to show that
  \[
  (\mu_r)_!( \nuMr \circ \dc_{M_r}\inverse \circ (\beta_{\mu_r}
  \inverse)_\sharp) \circ \phi_{\mu_r} \circ L_r(\sigma \inverse)\chk
  = L_r(\sigma) \circ (\mu_r)_!( \nuMr \circ \dc_{M_r}\inverse \circ
  (\beta_{\mu_r}\inverse )_\sharp )\circ \phi_{\mu_r}.
  \]
  In other words, we need to show that the diagram of local systems
  \begin{equation}
    \label{eq:lc}
    \xymatrix{ ((\mu_r)_! \BBQ_{M_r}) \chk [-2d] \ar[rrr]^{ L_r(\sigma
      \inverse) \chk} \ar[d]_{ (\mu_r)_!( (\beta_{\mu_r}\inverse
      )^\sharp) \circ \phi_{\mu_r}} &&& ((\mu_r)_!  \BBQ_{M_r}) \chk
    [-2d] \ar[d]^{ (\mu_r)_!((\beta_{\mu_r}\inverse )^\sharp) \circ
      \phi_{\mu_r})} \\ 
      (\mu_r)_! \BBQ_{M_r} \chk [-2d] \ar[d]_{ (\mu_r)_!(
        \dc_{M_r}\inverse)} &&& (\mu_r)_!  \BBQ_{M_r} \chk [-2d]
      \ar[d]^{ (\mu_r)_! (\dc_{M_r}\inverse)} \\
      (\mu_r)_! \BBD_{M_r} [-2d] \ar[d]_{ (\mu_r)_!( \nuMr)} &&&
      (\mu_r)_!  \BBD_{M_r} [-2d] \ar[d]^{(\mu_r)_! (\nuMr)} \\
      (\mu_r)_! \BBQ_{M_r} \ar[rrr]^{L_r(\sigma)}&&&
      (\mu_r)_!  \BBQ_{M_r} }  
  \end{equation}
  commutes.
  
  Using that $((\mu_r)_! \BBQ_{M_r}) \chk [-2d]$ is isomorphic to the
  dual local system, $((\mu_r)_! \BBQ_{M_r})^*$, and $(\mu_r)_!
  \BBQ_{M_r} \chk [-2d]$ is isomorphic to $(\mu_r)_!  \BBQ_{M_r}^*$,
  since $M_r$ and $N_r$ are rational homology manifolds, it is
  straightforward to show that for $x$ in $N_r$, the diagram obtained
  from diagram (\ref{eq:lc}) by taking the stalk at $x$ commutes. It
  follows that diagram (\ref{eq:lc}) commutes as desired.
\end{proof}

\subsection{}\label{4.7}
In this subsection we show that if $\Sigma \times \Sigma$ acts on
$\Ext_{N_0}^\bullet \left( \jN^* R\mu_!  \BBQ_{M}, \jN^!  R\mu_!
  \BBD_{M}\right)$ by
\[
(\sigma, \sigma')\cdot f= \jN^!L_Q'(\sigma') \circ f \circ
\jN^*L(\sigma \inverse) = \left(\jN^*L(\sigma \inverse)^\sharp \circ
\jN^!L_Q'(\sigma')_\sharp \right) (f),
\]
then $\can\inverse \circ (\nat_{\delta \jN} \circ \jN^!  \delta^!
(\lambda))_\sharp$ is $\Sigma\times \Sigma$-equivariant.

Suppose $\sigma$ and $\sigma'$ are in $\Sigma$ and $f$ is in
$\Ext_{N_0}^\bullet \left(\BBQ_{N_0}, \jN^!\delta^! ((R\mu_!
  \BBQ_{M})\chk \boxtimes R\mu_! \BBD_{M})\right)$. Then, setting $u=
L(\sigma \inverse)$, $v= L_Q'(\sigma')$ and using $\nat_{\delta
  \jN}= \nat_{\jN} \circ \jN^!( \nat_\delta)$, Corollary \ref{cor2.3},
the naturality of $\nat_{\jN}$, and Proposition \ref{prop2.3} we have:
\begin{align*}
  \can\inverse \circ \left( \nat_{\delta \jN} \circ \jN^!
    \delta^!(\lambda)\right)_\sharp ( (\sigma, \sigma')\cdot f) &=
  \can \inverse\left( \nat_{\delta \jN} \circ \jN^! \delta^!(\lambda
    \circ
    (u\chk \boxtimes v) ) \circ f\right)\\
  &= \can \inverse\left( \nat_{\jN} \circ \jN^!( \nat_\delta \circ
    \delta^!( \lambda \circ (u\chk \boxtimes v) )) \circ f \right) \\
  &=\can \inverse\left( \nat_{\jN} \circ \jN^!( u^\sharp \circ
    v_\sharp \circ \nat_\delta \circ \delta^!(\lambda) ) \circ f
  \right) \\
  &= \can \inverse\left( \jN^*(u)^\sharp \circ \jN^!(v)_\sharp \circ
    \nat_{\jN} \circ \jN^! ( \nat_\delta \circ \delta^!(\lambda))
    \circ f\right) \\
  &= \can \inverse\left( (\jN^*(u)^\sharp \circ
    \jN^!(v)_\sharp)_\sharp ( \nat_{\delta \jN} \circ \jN^!
    \delta^!(\lambda) \circ f) \right) \\
  &= \jN^*(u)^\sharp \circ \jN^!(v)_\sharp \circ \can\inverse
  (\nat_{\delta \jN} \circ \jN^!\delta^!(\lambda) \circ f) \\
  &= (\sigma, \sigma') \cdot \left( \can\inverse \circ (\nat_{\delta
      \jN} \circ \jN^! \delta^!(\lambda))_\sharp (f) \right) .
\end{align*}

\subsection{}\label{4.8}
To complete the proof of Theorem \ref{thm4.4} we need to show that 
\[
\xymatrix{\Ext_{N_0}^\bullet \left( \jN^* R\mu_! \BBQ_{M}, \jN^!
    R\mu_! \BBD_{M}\right) \ar[rr]^-{(a\inverse)^\sharp \circ b_\sharp}
  && \Ext_{N_0}^{\bullet+4n} \left( R(\mu_0)_!  \BBQ_{M_0}, R(\mu_0)_!
    \BBQ_{M_0}\right)}
\]
is $\Sigma\times \Sigma$-equivariant where $\Sigma \times \Sigma$ acts
on the domain as in \S4.7. For this, it is enough to show that
\[
b\circ \jN^!L_Q'(\sigma') \circ f\circ \jN^*L(\sigma \inverse) \circ
a\inverse = L_0(\sigma')\circ b\circ f\circ a\inverse \circ L_0(\sigma
\inverse)
\]
for $\sigma$ and $\sigma'$ in $\Sigma$ and $f$ in $\Ext_{N_0}^\bullet
\left( \jN^* R\mu_! \BBQ_{M}, \jN^!  R\mu_!  \BBD_{M}\right)$.

Recall from \S\ref{4.2} and \S\ref{4.4} that $a=R(\mu_0)_!(
\alpha_{\jM})\circ \bc^*$ and
\[
L_0(\sigma \inverse)= R(\mu_0)_!( \alpha_{\jM}) \circ \bc^*\circ
\jN^*L(\sigma \inverse) \circ (R(\mu_0)_!( \alpha_{\jM}) \circ \bc^*)
\inverse = a \circ \jN^*L(\sigma \inverse) \circ a\inverse.
\]

Therefore, to show that $(a\inverse)^\sharp \circ b_\sharp$ is
$\Sigma\times \Sigma$-equivariant, it is enough to show that $b\circ
\jN^!L_Q'(\sigma')= L_0(\sigma') \circ b$.

Recall from \S\ref{4.4} that $b=R(\mu_0)_!( \nu_{M_0} \circ \beta_{\jM}
\inverse )\circ \bc^!$, so we need to show that
\begin{multline*}
  R(\mu_0)_!( \nu_{M_0} \circ \beta_{\jM} \inverse )\circ \bc^! \circ
  \jN^!( R\mu_!(\nuM^Q) \inverse\circ L(\sigma) \circ R\mu_!(\nuM^Q))=
  \\
  R(\mu_0)_!( \alpha_{\jM}) \circ \bc^*\circ \jN^*L(\sigma) \circ
  (\bc^*)\inverse \circ R(\mu_0)_!( \alpha_{\jM}\inverse) \circ
  R(\mu_0)_!( \nu_{M_0} \circ \beta_{\jM} \inverse )\circ \bc^!
\end{multline*}
for $\sigma$ in $\Sigma$. 

Setting $\nuM= \nuM^Q$, and using the naturality of $\bc^!$, it
follows that it is enough to show that
\begin{multline*}
  \jN^!L(\sigma) \circ (\bc^!)\inverse \circ R(\mu_0)_! (\jM^!( \nuM)
  \circ \beta_{\jM} \circ \nu_{M_0} \inverse \circ \alpha_{\jM}) \circ
  \bc^* = \\
  (\bc^!)\inverse \circ R(\mu_0)_! (\jM^!( \nuM) \circ \beta_{\jM}
  \circ \nu_{M_0} \inverse \circ \alpha_{\jM}) \circ \bc^* \circ
  \jN^*L(\sigma).
\end{multline*}

Set $\tau= (\bc^!)\inverse \circ R(\mu_0)_! (\jM^!( \nuM) \circ
\beta_{\jM} \circ \nu_{M_0} \inverse \circ \alpha_{\jM}) \circ \bc^*$.
Then $\tau$ is an isomorphism in $\Dbc(N_0)$, $\tau\colon \jN^*R\mu_!
\BBQ_M \to \jN^! R\mu_!  \BBQ_M[2l]$, where $l=\codim_NN_0= \codim_
MM_0$, and we need to show that
\begin{equation}
  \label{eq:tau}
   \jN^! L(\sigma) \circ \tau = \tau \circ \jN^* L(\sigma).
\end{equation}

We prove the following proposition in the Appendix.

\begin{proposition}\label{tau=rho}
  There is a natural transformation, $\rho^{\jN}\colon \jN^*\to
  \jN^![2l]$, so that $\tau= \rho^{\jN}_{R\mu_!\BBQ_M}$.
\end{proposition}

Given the truth of the proposition, it follows from the naturality of
$\rho^{\jN}$ that $\jN^!(g)\circ \tau = \tau \circ \jN^*(g)$ for $g$
in $\End_{N}( R\mu_!  \BBQ_M)$ and so in particular (\ref{eq:tau})
holds for $\sigma$ in $\Sigma$. This completes the proof of Theorem
\ref{thm4.4}. 

\section{Generalized Steinberg Varieties}

\subsection{}\label{5.1}
In this section we apply the results of \S3 and \S4 to generalized
Steinberg varieties.

We start with the following incarnation of the basic commutative
diagram (\ref{eq:bd}) as in \cite{borhomacpherson:partial}:
\begin{equation} \label{eq:bdg}
  \xymatrix{ 
    \CNt \ar[rr]^-{\eta^\CP_0} \ar[d]_{\jNt}&& \CNtP
    \ar[rr]^-{\xi^\CP_0} \ar[d]_{\jNtP}&& \CN \ar[d]_{\jCN} \\ 
    \fgt \ar[rr]^-{\eta^\CP}&& \fgtP \ar[rr]^-{\xi^\CP}&& \fg\\
    \fgrst \ar[rr]^-{\etarsP} \ar[u]&& \fgrstP \ar[rr]^-{\xirsP} \ar[u]&&
    \fgrs \ar[u]}
\end{equation}
The notation is as follows:
\begin{itemize}
\item $G$ is a connected, reductive, complex algebraic group with Lie
  algebra $\fg$, $\CN$ is the cone of nilpotent elements in $\fg$, and
  $\fgrs$ is the open subvariety of regular semisimple elements in
  $\fg$.
\item $\CP$ is a conjugacy class of parabolic subgroups of $G$.
\item $\fgt= \{\, (x, B)\in \fg\times \CB \mid x\in \Lie(B) \,\}$
  where $\CB$ is the variety of Borel subgroups of $G$ and $\fgtP=
  \{\, (x, P)\in \fg\times \CP \mid x\in \Lie(P) \,\}$.
\item The maps $\eta^\CP_?$ are defined by $\eta^\CP_?(x,B)= (x,P)$
  where $P$ is the unique subgroup in $\CP$ that contains $B$.
\item The maps $\xi^\CP_?$ are projection on the first factor.
\item $\mu= \xi^\CP \circ \eta^\CP$ is the projection on the first
  factor.
\item $\fgrst= \mu\inverse( \fgrs)= \{\, (x, B)\in \fgrs\times \CB
  \mid x\in \Lie(B) \,\}$ and $\fgrstP= (\xi^\CP)\inverse (\fgrs)=
  \{\, (x, P)\in \fgrs\times \CP \mid x\in \Lie(P) \,\}$.
\item $\CNt= \mu\inverse(\CN)= \{\, (x, B)\in \CN\times \CB \mid x\in
  \Lie(B) \,\}$ and $\CNtP= (\xi^\CP)\inverse (\CN)= \{\, (x, P)\in
  \CN\times \CP \mid x\in \Lie(P) \,\}$.
\end{itemize}

In accordance with the notation above, we assume also that $\dim G=
\dim \fg= d$, $\dim \CN= 2n$, and $l= d-2n = \codim_{\fg} \CN$.

It is shown in \cite{borhomacpherson:partial} that diagram
(\ref{eq:bdg}) has properties D1 -- D7 of the basic commutative
diagram. For the convenience of the reader, we recall the group action
involved in properties D6 and D7.

Fix a maximal torus, $T$ and a Borel subgroup, $B$, of $G$ with
$T\subset B$. Define $\ft= \Lie(T)$ and $\ftrs= \ft\cap \fgrs$, so
$\ftrs$ is the set of regular semisimple elements in $\ft$. Let $W=
N_G(T)/T$ be the Weyl group of $(G,T)$. Then $W$ acts on $\ftrs\times
G/T$ on the right by $(t, gT) \cdot w= (\Ad(w\inverse) t, gwT)$ for
$w$ in $W$, $t$ in $\fgrst$ and $g$ in $G$. It is well-known and easy
to check that the rule $(t, gT)\mapsto (\Ad(g)t, gBg\inverse)$ defines
an isomorphism of varieties $\ftrs\times G/T\cong \fgrst$ and we use
this isomorphism to transport the $W$-action from $\ftrs\times G/T$ to
$\fgrst$. It is also well-known and easy to prove that the projection
on the first factor, from $\fgrst$ to $\fgrs$, is an orbit map for the
right $W$-action on $\fgrst$. Thus, diagram (\ref{eq:bdg}) has property
D6.

Next, let $P$ be the subgroup in $\CP$ with $B\subseteq P$ and set
$W_P=N_P(T)/T$, the Weyl group of $P$, so $W_P$ is a subgroup of $W$.
It is straightforward to check that $\eta^\CP|_{\fgrst}$ is an orbit
map for the action of $W_P$ on $\fgrst$. Thus, diagram (\ref{eq:bdg})
has property D7.

If $\CQ$ is a second conjugacy class of parabolic subgroups of $G$,
then the two variable version of diagram (\ref{eq:bdg}), as in \S3.5,
is the following:
\begin{equation} \label{eq:bdZ}
  \xymatrix{
    Z \ar[rr]^-{\eta} \ar[d]_{\jZ}&& X \ar[rr]^-{\xi} \ar[d]_{\jX}&&
    \CN \ar[d]_{\delta\circ \jCN} \\
    \fgt\times \fgt \ar[rr]^-{\eta^\CP\times \eta^\CQ}&& \fgtP\times
    \fgtQ \ar[rr]^-{\xi^\CP \times \xi^\CQ}&& \fg\times \fg\\
    \fgrst \times \fgrst \ar[rr] \ar[u]&& \fgrstP \times \fgrstQ
    \ar[rr] \ar[u]&& \fgrs \times \fgrs \ar[u]}
\end{equation}

Since $(\fgt \times \fgt) \times_{\fg \times \fg} \CN= \{\, ((x,B'),
(x,B'')) \mid x\in \Lie(B')\cap \Lie(B'') \,\}$, we may identify $Z$
with the Steinberg variety of $G$. Then $\jZ\colon Z\to \fgt \times
\fgt$ by $\jZ(x,B', B'')= ((x,B'), (x, B''))$.

Also, since $(\fgtP \times \fgtQ) \times_{\fg \times \fg} \CN = \{\,
((x,P'), (x,Q'))\mid x\in \Lie(P')\cap \Lie(Q') \,\}$, we may identify
$X$ with the generalized Steinberg variety $\XPQ$ from \S1. Then
$\jX\colon \XPQ \to \fgtP \times \fgtQ$ by $\jX(x,P', Q')= ((x,P'),
(x, Q'))$.

Applying Theorem \ref{cor3.5} we have our first main result.

\begin{theorem}\label{bmaverage}
  If $H_\bullet(Z)$ is given the $W\times W$-action induced from the
  $W\times W$-action on $(\murs\times \murs)_! \BBQ_{\fgrst \times
    \fgrst}$, then there is an isomorphism of vector spaces,
  $H_\bullet( \XPQ) \cong H_\bullet(Z) ^{W_P\times W_Q}$, so that the
  diagram
  \begin{equation*}
     \xymatrix{
       H_\bullet(Z) \ar[rr]^{\eta_*} \ar[dr]_\Av & & H_\bullet(\XPQ)
       \ar@{<->}[dl]^\simeq \\   
       & H_\bullet(Z) ^{W_P\times W_Q} &
     }
   \end{equation*}
  commutes.
\end{theorem}

\subsection{}\label{5.2}
Now we consider the special case of Theorem \ref{bmaverage} when
$\bullet= 2\dim Z= 4n$ as in \S4.1. Borho and MacPherson
\cite{borhomacpherson:weyl} have shown that the $\BBQ$-algebra
homomorphism, $\BBQ W \to \End_{\CN}( R(\mu_0)_! \BBQ_{\CNt})$ from
\S4.2 is an isomorphism.  Therefore, from (\ref{eq:hom}) we get the
result originally proved by Kazhdan and Lusztig
\cite{kazhdanlusztig:topological} and strengthened by Chriss and
Ginzburg \cite{chrissginzburg:representation}.

\begin{theorem}\label{weq}
  If $W\times W$ acts on $\BBQ W$ by $(w,w') \cdot x= w'xw\inverse$,
  then there are $W\times W$-equivariant isomorphisms
  \[
  \xymatrix{\BBQ W \ar[r]^-{\simeq}& \End_{\CN}( R(\mu_0)_!
    \BBQ_{\CNt}) \ar[r]^-{\simeq}& H_{4n}(Z)\op.}
  \]
\end{theorem}

Recall that $e_P$ denotes the primitive idempotent in $\BBQ W_P$
corresponding to the trivial representation of $W_P$.  Since $\left(
  \BBQ W \right) ^{W_P\times W_Q} = e_Q \BBQ We_P$, the next corollary
follows immediately from Theorems \ref{bmaverage} and \ref{weq}.

\begin{corollary}\label{cor:ze}
  The $W\times W$-equivariant isomorphism $\BBQ W \xrightarrow{\simeq}
  H_{4n}(Z)\op$ in Theorem \ref{weq} induces an isomorphism between the
  subspace $e_Q \BBQ We_P$ of $\BBQ W$ and $H_{4n}(\XPQ)$, the top
  Borel-Moore homology group of the generalized Steinberg variety,
  $\XPQ$:
  \[
  \xymatrix{\BBQ W \ar[r]^{\simeq} \ar[d]_{\Av} & H_{4n}(Z)\op
    \ar[d]^{\eta_*} \\
    e_Q \BBQ We_P \ar[r]_{\simeq} &H_{4n}(\XPQ) }
  \]
\end{corollary}

\subsection{}\label{5.3}

In this subsection, we use Corollary \ref{cor:ze} to compute the
action of a simple reflection in $W$ on $H_{4n}(Z)$. What we prove is
the analog for $H_{4n}(Z)$ of the ``easy'' part of Hotta's
transformations for the action of a simple reflection in the
cohomology of a Springer fibre. Our argument is inspired by Hotta's
argument in \cite{hotta:local}.

It is well-known that $W$ indexes the $G$-orbits on $\CB\times \CB$
and that if $Z_w$ denotes the preimage of the orbit indexed by $w$ in
$W$ under the projection of $Z$ onto $\CB \times \CB$ given by the
projection on the second and third factors, then the dimension of
$Z_w$ is $2n$ and the irreducible components of $Z$ are the closures
of the $Z_w$'s (see \cite{steinberg:desingularization}).  Thus, if
$[\Zwbar]$ denotes the canonical class of $\Zwbar$ in $H_{4n}(Z)$, it
follows that $\{\, [\Zwbar] \mid w\in W\,\}$ is a basis of
$H_{4n}(Z)$.

Recall that we have fixed a Borel subgroup, $B$, of $G$ containing
$T$.  The choice of $B$ determines a set of Coxeter generators of $W$
and hence a length function and a partial order, the Bruhat order, on
$W$.

For the time being we fix a simple reflection, $s$, in $W$ and let
$\CP_s$ denote the conjugacy class of minimal parabolic subgroups of
$G$ determined by $s$. Then $\CP_s$ and $\CB$ are conjugacy classes of
parabolic subgroups of $G$ and we may consider $\eta_*\colon H_{4n}(Z)
\to H_{4n}(X^{\CP_s \times \CB})$.

Let $P_s$ be the subgroup in $\CP_s$ that contains $B$. It is shown in
\cite[\S3]{douglassroehrle:geometry} that if $w$ is in $W$, then $\dim
\eta(Z_w) =\dim Z_w$ if and only if $w$ is minimal in its $(W_{P_s},
W_B)$-double coset. Since $W_{P_s} =\{\, 1,s\,\}$ and $W_B=
\{\,1\,\}$, it follows that $w$ is minimal in its double coset if and
only if $sw>w$ in the Bruhat order. Therefore, if $sw<w$ we have
$\eta_*( [\Zwbar]) =0$. It follows that $\dim \ker \eta_* \geq |W|/2$.

On the other hand, by Corollary \ref{cor:ze}, we may identify $\eta_*$
with the averaging map onto the set of $W_{P_s}\times W_B$-invariants
in $\BBQ W$. In this case, the averaging map from $\BBQ W$ to $(\BBQ
W)^{W_{P_s} \times W_B}$ is $x\mapsto \frac 12(x+xs)$ and so its
kernel is $\{\, x\in \BBQ W\mid xs=-x\,\}$ and has dimension equal
$|W|/2$.  Therefore, the kernel of $\eta_*$ is the subspace $\{\, c\in
H_{4n}(Z) \mid s\cdot c= -c\,\}$ and it has dimension equal $|W|/2$.
Since $\ker \eta_*$ contains the linearly independent set $\{\,
[\Zwbar]\mid sw<w \,\}$, it follows that $\{\, [\Zwbar]\mid sw<w \,\}$
is a basis of $\ker \eta_*$.  This proves the following theorem.

\begin{theorem}\label{hotta}
  If $s$ is a simple reflection in $W$, then $\{\, [\Zwbar]\mid sw<w
  \,\}$ is a basis of the subspace $\{\, c\in H_{4n}(Z) \mid s\cdot c=
  -c\,\}$ of $H_{4n}(Z)$. In particular, if $w$ is in $W$ and $s$ is a
  simple reflection, then $s\cdot [\Zwbar]= -[\Zwbar]$ if and only if
  $sw<w$ in the Bruhat order.
\end{theorem}

\subsection{}\label{5.4}

We now turn to computing the top Borel-Moore homology group of the
generalized Steinberg variety $\YPQ$. Recall that we have fixed
parabolic subgroups, $P$ in $\CP$ and $Q$ in $\CQ$, with $B\subseteq
P\cap Q$. Then
\[
\YPQ=\{\, (x, P', Q')\in \CN\times \CP\times \CQ \mid x\in
\Lie(U_{P'}) \cap \Lie(U_{Q'}) \,\} \subseteq \XPQ
\] 
and 
\[
\ZPQ= \eta\inverse(\YPQ).
\]
Thus, we have a cartesian square
\[
\xymatrix{\ZPQ \ar[r]^-{j} \ar[d]_{\etabar}&  Z \ar[d]^{\eta}\\
  \YPQ \ar[r]& \XPQ}
\]
where the horizontal arrows are inclusions and $\etabar$ is the
restriction of $\eta$ to $\ZPQ$.

It follows from the definitions that $\etabar$ is a fibre bundle with
smooth fibres isomorphic to $P/B \times Q/B$.

Define $\WPQ$ to be the set of maximal length $(W_P, W_Q)$-double
coset representatives in $W$, so $\WPQ$ indexes the $G$-orbits on
$\CP\times \CQ$.

It was shown in \cite[\S4]{douglassroehrle:geometry} that if $Y_w$
denotes the preimage of the orbit indexed by $w$ in $\WPQ$ under the
projection of $\YPQ$ onto $\CP \times \CQ$ given by the projection on
the second and third factors, then the dimension of $Y_w$ is $\dim
\CP+ \dim \CQ$ and the irreducible components of $\YPQ$ are the
closures of the $Y_w$'s.

It was also shown in \cite[\S4]{douglassroehrle:geometry} that $\{\,
\Zwbar\mid w\in \WPQ\,\}$ is the set of irreducible components of
$\ZPQ$. Clearly $\etabar(Z_w)\subseteq Y_w$ and so since $\etabar$ is
proper, $Z_w$ and $Y_w$ are irreducible, and the fibres $\etabar$ all
have the same dimension, it follows that $\etabar( \Zwbar)= \Ywbar$.

Since $\etabar$ is a fibre bundle with smooth fibres, if $f=\dim
P/B + \dim Q/B$, then there is an inverse image map in
Borel-Moore homology, $\etabar^*\colon H_\bullet(\YPQ) \to
H_{\bullet+2f}(\ZPQ)$ (see \cite[8.3.31]
{chrissginzburg:representation}).

It is straightforward to check that if $[\Ywbar]$ denotes the
canonical class of $\Ywbar$ in $H_{4n-2f}(\YPQ)$, then $\etabar^*(
[\Ywbar])$ is a multiple of $[\Zwbar]$ (see
\cite{fulton:intersection}). Since $\dim H_{4n-2f}(\YPQ)= \dim
H_{4n}(\ZPQ)$ it follows that $\etabar^*$ is injective.

Next, $\ZPQ$ is a closed subvariety of $Z$, so if $j$ denotes the
inclusion, there is a direct image map in Borel-Moore homology,
$j_*\colon H_\bullet(\ZPQ)\to H_\bullet(Z)$. It follows immediately
that $j_*( [\Zwbar])= [\Zwbar]$ for $w$ in $\WPQ$ and that $j_*$ is
injective. 

Combining the results in the last two paragraphs we have proven the
next proposition.

\begin{proposition}
  The mapping $\etabar^*\colon H_{4n-2f}(\YPQ) \to H_{4n}(\ZPQ)$ is an
  isomorphism of vector spaces and the mapping $j_*\colon H_{4n}(
  \ZPQ) \to H_{4n}(Z)$ is injective with image equal the span of $\{\,
  [\Zwbar] \mid w\in \WPQ\,\}$.
\end{proposition}

\subsection{}\label{5.5}

We identify the image of $H_{4n}(\ZPQ)$ with its image in $H_{4n}(Z)$.
Then $H_{4n}(\ZPQ)$ is the span of $\{\, [\Zwbar] \mid w\in \WPQ\,\}$
in $H_{4n}(Z)$ and $H_{4n}(\ZPQ) \cong H_{4n-2f}(\YPQ)$.  Define $\HPQ$
to be the subspace of $c$ in $H_{4n}(Z)$ with the property that
$s\cdot c=-c$ and $c\cdot t=-c$ for all simple reflections, $s$ in
$W_P$ and $t$ in $W_Q$. It follows from Theorem \ref{hotta} that
$H_{4n}(Z) \subseteq \HPQ$.

Recall that $\epsilon_P$ and $\epsilon_Q$ denote the primitive
idempotents in $W_P$ and $W_Q$ corresponding to the sign
representations of $W_P$ and $W_Q$ respectively. Then $\dim \epsilon_Q
\BBQ W \epsilon_P = |\WPQ|$ and $\epsilon_Q \BBQ W \epsilon_P$ is the
set of all $x$ in $\BBQ W$ with the property that $sx=-x$ and $xt=-x$
for all simple reflections, $s$ in $W_Q$, and $t$ in $W_P$. It follows
from Theorem \ref{weq} that under the isomorphism $\BBQ
W\xrightarrow{\simeq} H_{4n}(Z)\op$, the subspace $\HPQ$ is the image
of $\epsilon_Q \BBQ W \epsilon_P$. Therefore, $\dim \HPQ= | \WPQ|$ and
hence $\HPQ = H_{4n}(\ZPQ)$. This proves the following theorem.

\begin{theorem}
  The $W\times W$-equivariant isomorphism $\BBQ W \xrightarrow{\simeq}
  H_{4n}(Z)\op$ in Theorem \ref{weq} induces an isomorphism between the
  subspace $\epsilon_Q \BBQ W \epsilon_P$ of $\BBQ W$ and
  $H_{4n-2f}(\YPQ)$, the top Borel-Moore homology group of the
  generalized Steinberg variety, $\YPQ$,
  \[
  \xymatrix{ \epsilon_Q \BBQ W \epsilon_P \ar[r]^-\simeq \ar[d] &
    H_{4n-2f}(\YPQ) \ar[d]^{j_* \etabar^*} \\
    \BBQ W \ar[r]_{\simeq} & H_{4n}(Z)\op }
  \]
  where the left vertical arrow is inclusion.
\end{theorem}

\appendix
\section{}

\subsection{}\label{a.1}

In this appendix we change notation slightly from \S\ref{2.4}. For a
morphism, $\xi\colon X\to Y$, of complex, algebraic varieties, the
units of the adjoint pairs $(\xi^*, \xi_*)$ and $(\xi^!, \xi_!)$ are
denoted by $\eta_\xi^*$ and $\eta_\xi^!$ respectively. Similarly, the
counits are denoted by $\epsilon_\xi^*$ and $\epsilon_\xi^!$
respectively.

Suppose $\xi\colon X\to Y$ is a morphism between complex, algebraic
varieties that are rational homology manifolds. For a fixed choice of
isomorphisms $\nu_X\colon \BBD_X \xrightarrow{\simeq} \BBQ_X[2\dim X]$
and $\nu_Y\colon \BBD_Y \xrightarrow{\simeq} \BBQ_Y[2\dim Y]$ there is
a natural transformation $\rho^\xi\colon \xi^*\to \xi^![2l]$ where $l=
\dim Y- \dim X$. For a complex $A$ in $\Dbc(Y)$, $\rho^\xi=\rho^\xi_A$
is defined to be the composition
\begin{multline*}
  \xymatrix{\xi^*A \ar[r]^-{m_1\inverse} & \BBQ_X\otimes \xi^*A
    \ar[rr]^-{\omega_\xi \otimes \id} && \xi^!\BBQ_X\otimes \xi^*A[2l]
    \ar[r]^-{\eta_\xi^!} &\xi^!R\xi_!( \xi^!\BBQ_X\otimes
    \xi^*A) [2l] \ar[rr]^-{\xi^!(\pr_\xi\inverse)} && {}} \\
  \xymatrix{\xi^!( R\xi_!\xi^!\BBQ_X\otimes A) [2l]
    \ar[rr]^-{\xi^!(\epsilon_\xi^! \otimes \id)} && \xi^!( \BBQ_X
    \otimes A)[2l] \ar[rr]^-{\xi^!(m_1)} &&\xi^!A [2l]}
\end{multline*}
where the notation is as follows:
\begin{itemize}
\item $m_1\colon \BBQ_X \otimes B \xrightarrow{\simeq} B$ is the natural
  isomorphism for $B$ in $D(X)$.
\item $\omega_\xi= \xi^!(\nu_Y) \circ \beta_\xi \circ \nu_X\inverse
  \colon \BBQ_X \xrightarrow{\simeq} \xi^! \BBQ_Y[2l]$, so $\omega_\xi$
  is an isomorphism in $\Dbc(X)$.
\item $\eta_\xi^!$ and $\epsilon_\xi^!$ are as above.
\item For $B$ in $\Db(X)$ and $C$ in $\Db(Y)$, $\pr_\xi\colon R\xi_!
  B\otimes C \xrightarrow{\simeq} R\xi_!(B\otimes \xi^*C)$ is the
  projection isomorphism.
\end{itemize}
Notice that $\rho^\xi$ is a natural transformation since each map in
the definition of $\rho^\xi_A$ is natural in $A$.

Now consider a cartesian square
\begin{equation}
  \label{eq:cd}
  \xymatrix{M_0 \ar[r]^{\mu_0} \ar[d]_{\jM} & N_0 \ar[d]^{\jN} \\
  M \ar[r]_\mu & N}  
\end{equation}
satisfying the following conditions:
\begin{itemize}
\item[C1] The spaces are all complex, algebraic varieties that are
  rational homology manifolds.
\item[C2] The maps are all proper morphisms.
\item[C3] $\jM$ and $\jN$ are closed embeddings.
\item[C4] $\dim M_0= \dim N_0= 2n$, $\dim M=\dim N=d$, and $l=d-2n$. 
\end{itemize}

For a cartesian square as in (\ref{eq:cd}), we have base change
isomorphisms
\[
\xymatrix{\bc^*\colon \jN^* R\mu_! \ar[r]^-{\simeq}& R(\mu_0)_!
  \jM^*} \qquad \text{and} \qquad \xymatrix{\bc^!\colon \jN^! R\mu_!
  \ar[r]^-{\simeq}& R(\mu_0)_!  \jM^!}
\]
defined as in \S\ref{2.5}.

We prove the following lemmas in the next two subsections.

\begin{lemma}\label{lema1}
  If $X$ and $Y$ are complex, algebraic varieties that are rational
  homology manifolds and $\xi\colon X\to Y$ is a proper morphism, then
  $\rho^\xi_{\BBQ_Y} = \omega_\xi \circ \alpha_\xi$.
\end{lemma}

\begin{lemma}\label{lema2}
  If $\nu_{M_0}$ is chosen appropriately, then in the cartesian square
  (\ref{eq:cd}) the morphisms $\rho^{\jM}_{\BBQ_M}$ and $\rho^{\jN}
  _{R\mu_!\BBQ_M}$ are related by
  \[
  \bc^!\circ \rho^{\jN} _{R\mu_!\BBQ_M} = R(\mu_0)_! ( \rho^{\jM}
  _{\BBQ_M} ) \circ \bc^*.
  \]
\end{lemma} 

Recall from \S\ref{4.8} that in the setting of (\ref{eq:cd}) we have
$\tau\colon \jN^*R\mu_! \BBQ_M \xrightarrow{\simeq} \jN^! R\mu_!
\BBQ_M[2l]$, by
\[
\tau= (\bc^!)\inverse \circ R(\mu_0)_! (\jM^!( \nuM) \circ \beta_{\jM}
\circ \nu_{M_0} \inverse \circ \alpha_{\jM}) \circ \bc^*
\]
where $\nuM\colon \BBD_M \xrightarrow{\simeq} \BBQ_M[2d]$ and
$\nu_{M_0}\colon \BBD_{M_0} \xrightarrow{\simeq} \BBQ_{M_0}[4n]$ are
isomorphisms in $\Dbc(M)$ and $\Dbc(M_0)$ respectively.

Assuming Lemmas \ref{lema1} and \ref{lema2} have been proved we
have: 
\begin{align*}
  \rho^{\jN}_{R\mu_!\BBQ_M} &= (\bc^!)\inverse \circ R(\mu_0)_! (
  \rho^{\jM}_{\BBQ_M} ) \circ \bc^*\\
  &= (\bc^!)\inverse \circ R(\mu_0)_! ( \omega_{\jM} \circ
  \alpha_{\jM} ) \circ \bc^* \\
  &= (\bc^!)\inverse \circ R(\mu_0)_! ( \jM^!(\nuM) \circ \beta_{\jM}
  \circ \nu_{M_0}\inverse \circ
  \alpha_{\jM} ) \circ \bc^* \\
  &= \tau
\end{align*}
This proves Proposition \ref{tau=rho}.

\subsection{}\label{a.2}
In this subsection we prove Lemma \ref{lema1}.  Before doing so, we
need some preliminary results.

If $A$ is in $\Dbc(X)$ we denote the canonical isomorphisms
$\BBQ_X\otimes A \xrightarrow{\simeq} A$ and $A\otimes \BBQ_X
\xrightarrow{\simeq} A$ by $m_1$ and $m_2$ respectively. When
$A=\BBQ_X$ we set $m=m_1=m_2$.

The proof of the next lemma is a straightforward computation using
stalks and is omitted.

\begin{lemma}\label{lema.2.1}
  Suppose $A$ and $B$ are in $\Dbc(X)$ and $p_A\colon A\to \BBQ_X$ and
  $p_B\colon B\to \BBQ_X$ are two morphisms in $\Dbc(X)$. Then the
  diagrams
  \[
  \begin{gathered}
    \xymatrix{ R\xi_!A \otimes \BBQ_Y \ar[r]^-{\pr_\xi} \ar[d] _{m_2}&
      R\xi_!( A \otimes \xi^* \BBQ_Y )
      \ar[d]^{R\xi_!(\id\otimes \alpha_\xi)} \\
      R\xi_!A & R\xi_!( A \otimes \BBQ_X ) \ar[l] _{R\xi_!(m_2)} }
  \end{gathered}
  \quad \text{and} \quad
  \begin{gathered}
    \xymatrix{A\otimes B \ar[r]^{\id\otimes p_B} \ar[d]_{ p_A\otimes
        \id} & A\otimes \BBQ_X \ar[r]^{m_2} \ar[d]^{p_A\otimes \id} &A
      \ar[dd] ^{p_A} \\
      \BBQ_X\otimes B\ar[r]_{\id\otimes p_B} \ar[d]_{m_1} &\BBQ_X
      \otimes \BBQ_X \ar[rd]_m & \\
      B\ar[rr]_{p_B} &&\BBQ_X}
  \end{gathered}
  \]
  commute.
\end{lemma}

Let $\natt_\xi\colon \xi^*(A\otimes B) \xrightarrow{\simeq} \xi^*A
\otimes \xi^*B$ denote the canonical isomorphism in $\Db(X)$.

\begin{lemma}\label{lema.2.2}
  The diagram
  \[
  \xymatrix{R\xi_! \xi^! \BBQ_Y \otimes \BBQ_Y \ar[r]^-{\pr_\xi}
    \ar[dd] _{\epsilon_\xi^! \otimes \id}& R\xi_!( \xi^! \BBQ_Y
    \otimes \xi^* \BBQ_Y ) \ar[d]^{R\xi_!(\id\otimes \alpha_\xi)} \\
    & R\xi_!( \xi^! \BBQ_Y \otimes \BBQ_X ) \ar[d]
    ^{\Phi_\xi\inverse (m_2)} \\
    \BBQ_Y \otimes \BBQ_Y \ar[r]_-{m} &\BBQ_Y }
  \]
  commutes.
\end{lemma}

\begin{proof}
  Using the definition of $\Phi_\xi\inverse$ we have 
  \[
  \Phi_\xi\inverse(m_2) \circ R\xi_!(\id\otimes \alpha_\xi) \circ
  \pr_\xi = \epsilon_\xi^! \circ R\xi_!(m_2) \circ R\xi_!(\id\otimes
  \alpha_\xi) \circ \pr_\xi .
  \]
  Also, using the naturality of $\epsilon_\xi^!$ we have $m\circ
  (\epsilon_\xi^! \otimes \id)= \epsilon_\xi^! \circ m_2$, so it is
  enough to show that
  \begin{equation}
    \label{eq:a.2.1}
    m_2= R\xi_!(m_2) \circ R\xi_!(\id\otimes \alpha_\xi) \circ \pr_\xi.    
  \end{equation}
  
  Since $\xi$ is proper, we have $\pr_\xi= \Psi_\xi( (\epsilon_\xi^*
  \otimes \id) \circ \natt_\xi )$.
  
  The proof of (\ref{eq:a.2.1}) is a straightforward computation using
  the formula for $\pr_\xi$ and Lemma \ref{lema.2.2}. We omit the
  details.
\end{proof}

We can now complete the proof of Lemma \ref{lema1}. Recall that 
\[
\rho^\xi_{\BBQ_Y}= \xi^!\left( m_1\circ \epsilon_\xi^! \circ
  \pr_\xi\inverse \right) \circ \eta_\xi^! \circ (\omega_\xi \otimes
\id) \circ m_1\inverse = \Phi_\xi(m_1\circ \epsilon_\xi^! \circ
\pr_\xi\inverse) \circ (\omega_\xi \otimes \id) \circ m_1\inverse,
\]
so to prove the lemma, we need to show that
\[
\Phi_\xi(m_1\circ \epsilon_\xi^! \circ \pr_\xi\inverse) = \omega_\xi
\circ \alpha_\xi \circ m_1 \circ (\omega_\xi \inverse \otimes \id).
\]

Taking $A= \xi^!\BBQ_Y[2l]$, $p_A= \omega_\xi\inverse$, $B= \xi^*
\BBQ_Y$, and $p_B= \alpha_\xi$ in Lemma \ref{lema.2.1} we get
\[
\alpha_\xi \circ m_1\circ \circ (\omega_\xi \inverse \otimes \id)=
\omega_\xi\inverse \circ m_2 \circ (\id\otimes \alpha_\xi),
\]
so it is enough to show that $\Phi_\xi(m_1\circ \epsilon_\xi^! \circ
\pr_\xi\inverse) = m_2 \circ (\id\otimes \alpha_\xi)$. This last
equality follows immediately from Lemma \ref{lema.2.2}. This completes
the proof of Lemma \ref{lema1}.

\subsection{}\label{a.3}
In this subsection we prove Lemma \ref{lema2}. 

The proof is accomplished by showing that the diagrams (\ref{d:a1})
and (\ref{d:a2}) below are commutative. Then juxtaposing these
diagrams and tracing around the outside gives the desired result.

It is easy to see that any unlabeled regions of diagrams (\ref{d:a1})
and (\ref{d:a2}) commute. The commutativity of the labeled regions is
shown in the corresponding statements below.

To make the diagrams as clear as possible, we need to simplify the
notation. First, for a morphism, $\xi\colon X\to Y$, we denote the
derived functors $R\xi_*$ and $R\xi_!$ simply by $\xi_*$ and $\xi_!$
respectively. Second, we denote $\jN$ simply by $j$.  Third, we label
the maps in the diagrams using only the core maps or natural
transformations involved. For example, we write $\alpha_{\mu_0}$
instead of $(\mu_0)_!(\alpha_{\mu_0} \otimes \id)$ and $\bc^*$ instead
of $j^!j_!(\id\otimes \bc^*)$.

If $\xi\colon X\to Y$, and $A$ and $B$ are complexes, then
\[
\xymatrix{\pr^1\colon \xi_!A\otimes B \ar[r]^-{\simeq}& \xi_!(A\otimes
  \xi^*B)} \quad \text{and} \quad \xymatrix{\pr^2\colon A\otimes
  \xi_!B \ar[r]^-{\simeq}& \xi_!(\xi^*A\otimes B).}
\]
With this notation we have $\rho^\xi= \xi^!(m_1 \circ (\epsilon_\xi^!
\otimes \id) \circ (\pr^1_\xi) \inverse) \circ \eta_\xi^!  \circ
(\omega_\xi \otimes \id) \circ m_1\inverse$.

Notice that if $\xi$ is proper, then $\pr_\xi^1= \Psi_\xi(
(\epsilon_\xi^* \otimes \id) \circ \natt_\xi )$ and $\pr_\xi^2=
\Psi_\xi( (\id \otimes \epsilon_\xi^*) \circ \natt_\xi )$.

For a cartesian square as in (\ref{eq:cd}) we have a base change
isomorphism, 
\[
\xymatrix{\bct\colon \mu^* R(\jN)_!  \ar[r]^-{\simeq}& R(\jM)_!
  \mu_0^*,}
\]
defined by $\bct= \Psi_{\jM}( \mu_0^* (\epsilon_{\jN}^*))$.  Define
$\sigma\colon \mu_0^* \jN^!\to \jM^! \mu^*$ by $\sigma =\Phi_{\jM}(
\mu^*(\epsilon_{\jN}^!) \circ \bct \inverse)$. Then $\sigma$ is a
natural transformation.

\vfill\eject

\begin{equation}  \label{d:a1}
\begin{sideways}
\xymatrix{ 
j^*\mu_!\BBQ_M \ar[dddddddddd]^{\bc^*} 
&&&\BBQ_{N_0} \otimes j^* \mu_! \BBQ_M \ar[r]^{\omega_j} \ar[dd]^{\bc^*}
\ar[lll]_{m_1}
&j^! \BBQ_N \otimes j^*\mu_! \BBQ_M \ar[r]^{\eta_j^!} \ar[dd]^{\bc^*}
&j^!j_!(j^!\BBQ_N \otimes j^*\mu_!\BBQ_M) \ar[dd]^{\bc^*}  \\    
&&&&&\\
&&&\BBQ_{N_0} \otimes (\mu_0)_! \jM^* \BBQ_M \ar[r]^{\omega_j}
\ar[dd]^{\pr^2}  
&j^! \BBQ_N \otimes (\mu_0)_! \jM^* \BBQ_M \ar[r]^{\eta_j^!}
\ar[dd]^{\pr^2} &j^!j_!(j^!\BBQ_N \otimes (\mu_0)_! \jM^*\BBQ_M)
\ar[dd]^{\pr^2} \\  
&&&&&\\
&&&(\mu_0)_!( \mu_0^*\BBQ_{N_0} \otimes \jM^* \BBQ_M) \ar[r]^-{\omega_j}
\ar[dddddd]^{\alpha_{\mu_0}} 
&(\mu_0)_!((\mu_0)^* j^! \BBQ_N \otimes \jM^* \BBQ_M) \ar[r]^{\eta_j^!}
\ar[dd]^{\sigma} &
j^!j_!(\mu_0)_!(\mu_0^*j^!\BBQ_N \otimes  \jM^*\BBQ_M)
\ar[dd]^{\sigma} \\ 
&&&&&\\
&&&&(\mu_0)_!(\jM^! \mu^* \BBQ_N \otimes \jM^* \BBQ_M) \ar[r]^{\eta_j^!}
\ar[dddd]^{\alpha_\mu}
&j^!j_!(\mu_0)_!(\jM^! \mu^*\BBQ_N \otimes  \jM^*\BBQ_M)
\ar[dd]^{\alpha_\mu} \\  
&&&&&\\
&&\ar@{}[ur]^<{\text{(\ref{lema.2.1})}_1}&\ar@{}[r]^
{\text{(\ref{eq:a.3.2})}}& &j^!j_!(\mu_0)_! (\jM^! \BBQ_M \otimes
\jM^*\BBQ_M) \ar[dd]^{\bc^!}\\ 
&&&&&\\
(\mu_0)_! \jM^* \BBQ_M \ar@{<-} [uuuuuuuurrr]+L^{m_1}
&&&(\mu_0)_!(\BBQ_{M_0} \otimes \jM^* \BBQ_M) \ar[r]^{\omega_j}
\ar[lll]_-{m_1} 
&(\mu_0)_!(\jM^! \BBQ_M \otimes \jM^* \BBQ_M) \ar[r]^-{\eta_{\jM}^!}
\ar[ruu]^-{\eta_j^!} \ar@{}[ur]|(.8){\text{(\ref{eq:a.3.3})}}
&(\mu_0)_! \jM^!(\jM)_!(\jM^! \BBQ_M \otimes  \jM^*\BBQ_M)
}
\end{sideways}
\end{equation}

\vfill\eject

\begin{equation}  \label{d:a2}
\begin{sideways}
\xymatrix{ 
j^!j_!(j^!\BBQ_N \otimes j^*\mu_!\BBQ_M) \ar[dd]^{\bc^*} &
j^!(j_!j^!\BBQ_N \otimes \mu_!\BBQ_M) \ar[dd]^{\pr^2}
\ar[l]_{\pr^1} \ar[r]^{\epsilon_j^!} &j^!( \BBQ_N \otimes \mu_! \BBQ_M)
\ar[dd]^{\pr^2} \ar[rr]^-{m_1} &&j^!\mu_! \BBQ_M
\ar[dddddddddd]^{\bc^!} \ar@{<-} [ddddddddll]+R_{m_1} \\   
&&&&\\
j^!j_!(j^!\BBQ_N \otimes (\mu_0)_! \jM^*\BBQ_M) \ar[dd]^{\pr^2}
\ar@{}[uur]|{\text{(\ref{lema.3.1})}} & j^! \mu_!(\mu^*j_!j^!\BBQ_N
\otimes \BBQ_M) \ar[dd]^{\bct} \ar[r]^{\epsilon_j^!} &j^!\mu_!(
\mu^*\BBQ_N  \otimes \BBQ_M) \ar[dddd]^{=} \ar@{}
[urr]^{\text{(\ref{lema.2.1})}_2}   
  &&\\  
&&&&\\
j^!j_!(\mu_0)_!(\mu_0^*j^!\BBQ_N \otimes  \jM^*\BBQ_M)
\ar[dd]^{\sigma} & j^! \mu_!((\jM)_! \mu^*j^!\BBQ_N \otimes \BBQ_M)
  \ar[dd]^{\sigma} \ar[l]_{\pr^1} \ar@{}[uur]|
  {\text{(\ref{eq:a.3.1})}} && \\   
&&&&\\
j^!j_!(\mu_0)_!(\jM^! \mu^*\BBQ_N \otimes  \jM^*\BBQ_M)
\ar[dd]^{\alpha_\mu} & j^! \mu_!((\jM)_!
\jM^! \mu^*\BBQ_N \otimes \BBQ_M) \ar[dd]^{\alpha_\mu} \ar[l]_{\pr^1}
\ar[r]^-{\epsilon_{\jM}^!} &j^!\mu_!( \mu^*\BBQ_N \otimes \BBQ_M)
\ar[dd]^{\alpha_\mu} &&\\  
&&&&\\
j^!j_!(\mu_0)_!(\jM^! \BBQ_M \otimes  \jM^*\BBQ_M) \ar[dd]^{\bc^!} &
j^! \mu_!((\jM)_! \jM^! \BBQ_M \otimes 
\BBQ_M) \ar[dd]^{\bc^!} \ar[l]_{\pr^1} \ar[r]^-{\epsilon_{\jM}^!}
&j^!\mu_!( \BBQ_M \otimes \BBQ_M) \ar[dd]^{\bc^!} && \\
&&&&\\
(\mu_0)_! \jM^!(\jM)_!(\jM^! \BBQ_M \otimes  \jM^*\BBQ_M)& (\mu_0)_!
\jM^! ((\jM)_! \jM^! \BBQ_M 
\otimes \BBQ_M) \ar[l]_-{\pr^1} \ar[r]^-{\epsilon_{\jM}^!}  &(\mu_0)_!
\jM^!( \BBQ_M \otimes \BBQ_M) \ar[rr]^-{m_1} && (\mu_0)_! \jM^!\BBQ_M 
}
\end{sideways}
\end{equation}

\vfill\eject

The regions labeled $(\text{\ref{lema.2.1}})_1$ and
$(\text{\ref{lema.2.1}})_2$ commute using the analog of the first
rectangle in Lemma \ref{lema.2.1} with $m_1$ instead of $m_2$, taking
$A= \jM^!\BBQ_M$ and $A= \BBQ_M$ respectively.

\begin{lemma}\label{lema.3.2}
  The mapping $\sigma_{\BBQ_N} \colon \mu_0^* \jN^! \BBQ_N \to \jM^!
    \mu^* \BBQ_N$ is an isomorphism in $\Dbc(M_0)$.
\end{lemma}

\begin{proof}
  We have
  \[
  \sigma= \Phi_{\jM}( \mu^*(\epsilon_{\jN}^! \circ \bct \inverse) =
  \jM^! \mu^*( \epsilon_{\jN}^!) \circ \jM^!(\bct \inverse) \circ
  \eta_{\jM}^!.
  \]
  Since $\jM$ is a closed embedding, $\eta_{\jM}^!$ is an isomorphism,
  so it is enough to show that $\jM^! \mu^*( \epsilon_{\jN}^!) \colon
  \jM^!\mu^* (\jN)_! \jN^! \BBQ_N \longrightarrow \jM^!\mu^* \BBQ_N$
  is an isomorphism. Since $M$ and $N$ are purely $d$-dimensional,
  rational homology manifolds and $M_0$ and $N_0$ are purely
  $2n$-dimensional, rational homology manifolds, it follows that
  $\jM^!\mu^* (\jN)_!  \jN^! \BBQ_N$ and $\jM^!\mu^* \BBQ_N$ are both
  isomorphic to $\BBQ_{M_0} [-2l]$. It follows that $\jM^! \mu^*(
  \epsilon_{\jN}^!)$ is an isomorphism.
\end{proof}

Since $\sigma_{\BBQ_N}$ is an isomorphism, the composition
\[
\xymatrix{\beta_{\jM}\inverse \circ \jM^!( \nu_M\inverse \circ
  \alpha_\mu) \circ \sigma \circ \mu_0^*( \omega_{\jN} \inverse) \circ
  \alpha_\mu \inverse \colon \BBQ_{M_0} \ar[r]& \BBD_{M_0}[-4n]}
\]
is an isomorphism, so we may choose $\nu_{M_0}$ so
that 
\[
\nu_{M_0} \inverse = \beta_{\jM}\inverse \circ \jM^!( \nu_M\inverse
\circ \alpha_\mu) \circ \sigma_{\BBQ_N} \circ \mu_0^*( \omega_{\jN}
\inverse) \circ \alpha_\mu \inverse.
\]
It then follows that
\begin{equation}
  \label{eq:a.3.2}
  \omega_{\jM} \circ \alpha_{\mu_0} = \jM^!(\alpha_\mu) \circ
  \sigma_{\BBQ_N} \circ  \mu_0^*( \omega_{\jN}). 
\end{equation}

As in \S\ref{2.5}, $(\bc^!)\inverse \colon (\mu_0)_! \jM^!\to
\jN^!\mu_!$ by $(\bc^!)\inverse = \Phi_{\jN}(\mu_!(\epsilon_{\jM}^!))$.
Thus, using (\ref{eq:ad2}) we have
\begin{align*}
  (\bc^!)\inverse \circ (\mu_0)_!(\eta_{\jM}^!) &=
  \Phi_{\jN}(\mu_!(\epsilon_{\jM}^!)) \circ (\mu_0)_!(\eta_{\jM}^!)\\
  &= \Phi_{\jN}(\mu_!(\epsilon_{\jM}^!) \circ
  j_!(\mu_0)_!(\eta_{\jM}^!))\\
  &= \Phi_{\jN}(\mu_!(\epsilon_{\jM}^! \circ {\jM}_!(\eta_{\jM}^!)))\\
  &= \Phi_{\jN}(\mu_!(\Phi_{\jM}\inverse (\eta_{\jM}^!)))\\
  &= \Phi_{\jN}(\id)\\
  &= \eta_{\jN}^!.
\end{align*}
Therefore,
\begin{equation}
  \label{eq:a.3.3}
  \bc^! \circ \eta_{\jN}^!= (\mu_0)_!(\eta_{\jM}^!).
\end{equation}

\begin{lemma}\label{lema.3.1}
  The diagram
  \[
  \xymatrix{R(\jN)_!\jN^!A \otimes R\mu_!B \ar[rr]^{\pr_{\jN}^1} \ar[d]
    _{\pr_\mu^2} && R(\jN)_! \left( \jN^!A \otimes \jN^* R\mu_!B
    \right) \ar[d] ^{R(\jN)_!(\id \otimes \bc^*)} \\
    R\mu_! \left( \mu^* R(\jN)_! \jN^! A\otimes B \right) \ar[d]
    _{R\mu_!( \bct \otimes \id)} && R(\jN)_!\left( \jN^!A\otimes
      R(\mu_0)_! \jM^*B \right) \ar[d] ^{ R(\jN)_!( \pr_{\mu_0}^2)} \\
    R\mu_! \left( R(\jM)_! \mu_0^* \jN^! A\otimes B \right) \ar[rr]
    _{R\mu_!( \pr_{\jM}^!)} && R(\jN\mu_0)_!\left( \mu_0^*
      \jN^!A\otimes \jM^*B \right) }
  \]
  commutes for $A$ in $\Db(N)$ and $B$ in $\Db(M)$.
\end{lemma}

\begin{proof}
  First, using the formulas for $\pr_{\jN}^1$ and $\pr_\mu^2$ and the
  analogs of equation (\ref{eq:ad2}) for $\Psi_{\jN}$ and $\Psi_\mu$
  we see that it is enough to show that
  \begin{equation*}
    \Psi_\mu\left( \pr_{\jM}^1 \circ (\bct\otimes \epsilon_\mu^*)
      \circ \natt_\mu \right) =  \Psi_{\jN} \left( \pr_{\mu_0}^2 \circ
      (\epsilon_{\jN}^*\otimes \bc^* ) \circ \natt_{\jN} \right) .
  \end{equation*}

  Next, using the formulas for $\pr_{\jM}^1$ and $\pr_{\mu_0}^2$ and the
  analogs of equation (\ref{eq:ad2}) for $\Psi_{\jM}$ and
  $\Psi_{\mu_0}$ we see that it is enough to show that
  \begin{multline*}
    \Psi_\mu\Psi_{\jM} \left( (\epsilon_{\jM}^* \otimes \id) \circ
      \natt_{\jM} \circ \jM^*( (\bct\otimes \epsilon_\mu^*)
      \circ \natt_\mu )\right)\\
    = \Psi_{\jN} \Psi_{\mu_0} \left( (\id \otimes \epsilon_{\mu_0}^*)
      \circ \natt_{\mu_0} \circ \mu_0^*( (\epsilon_{\jN}^*\otimes \bc^*
      ) \circ \natt_{\jN}) \right) .
  \end{multline*}
  
  Now using that $\Psi_f \Psi_{g}= \Psi_{fg}$ and the naturality of
  $\natt_{\jM}$ and $\natt_{\mu_0}$, we see that it is enough to show
  that
  \begin{multline*}
    (\epsilon_{\jM}^* \otimes \id) \circ (\jM^* (\bct)\otimes \jM^*
    (\epsilon_\mu^*)) \circ \natt_{\jM} \circ \jM^*(\natt_\mu) \\=
    (\id \otimes \epsilon_{\mu_0}^*) \circ (\mu_0^* (\epsilon_{\jN}^*)
    \otimes \mu_0^* (\bc^*) ) \circ \natt_{\mu_0} \circ \mu_0^*
    (\natt_{\jN}) .
  \end{multline*}
  
  Since $\natt_{g} \circ g^*(\natt_f)= \nat_{fg}$, we only need to show
  that
  \[
  \epsilon_{\jM}^* \circ \jM^* (\bct)= \mu_0^* (\epsilon_{\jN}^*)
  \quad \text{and}\quad \jM^* (\epsilon_\mu^*) =\epsilon_{\mu_0}^*
  \circ \mu_0^* (\bc^*)
  \]
  which is the same as
  \[
  \Psi_{\jM} \inverse(\bct)= \mu_0^* (\epsilon_{\jN}^*) \quad
  \text{and}\quad \jM^* (\epsilon_\mu^*) = \Psi_{\mu_0}
  \inverse(\bc^*).
  \]
  
  These last two equations follow immediately from the definitions
  $\bct= \Psi_{\jM}( \mu_0^* (\epsilon_{\jN}^*))$ and $\bc^*=
  \Psi_{\mu_0}( \jM^* (\epsilon_\mu^*))$ above.
\end{proof}

Since $\Phi_{\jM}(\xi)= \jM^!(\xi) \circ \eta_{\jM}^!$, using the
naturality of $\epsilon_{\jM}^!$ and the fact that $\epsilon_{\jM}^!
\circ \eta_{\jM}^!= \id$, we have
\begin{equation}
  \label{eq:a.3.1}
  \begin{split} 
    \epsilon_{\jM}^! \circ (\jM)_!( \sigma) \circ \bct
    &=\epsilon_{\jM}^! \circ (\jM)_! \left(
      \jM^!(\mu^*(\epsilon_{\jN})
      \circ \bct\inverse) \right) \circ \eta_{\jM}^! \circ \bct \\
    &=\mu^*(\epsilon_{\jN}) \circ \bct\inverse \circ \epsilon_{\jM}^!
    \circ
    \eta_{\jM}^! \circ \bct \\
    &=\mu^*(\epsilon_{\jN}).
  \end{split}
\end{equation}



\begin{thebibliography}{10}

\bibitem{borel:intersection}
A.~Borel (ed.), \emph{Intersection cohomology ({B}ern, 1983)}, Progr. Math.,
  vol.~50, Birkh{\"a}user, Boston, 1984.

\bibitem{borhomacpherson:weyl}
W.~Borho and R.~MacPherson, \emph{Repr{\'e}sentations des groupes de {W}eyl et
  homologie d'intersection pour les vari{\'e}t{\'e}s nilpotentes}, C. R. Acad.
  Sci. Paris \textbf{292} (1981), no.~15, 707--710.

\bibitem{borhomacpherson:partial}
\bysame, \emph{Partial resolutions of nilpotent varieties. {A}nalysis and
  topology on singular spaces, {II}, {III} ({L}uminy, 1981)}, Ast\'erisque
  \textbf{101} (1983), 23--74.

\bibitem{chrissginzburg:representation}
N.~Chriss and V.~Ginzburg, \emph{Representation theory and complex geometry},
  Birkh{\"{a}}user, Boston, 1997.

\bibitem{douglassroehrle:geometry}
J.~M. Douglass and G.~R{\"o}hrle, \emph{The geometry of generalized {S}teinberg
  varieties}, Adv. Math. \textbf{187} (2004), no.~2, 396--416.

\bibitem{fulton:intersection}
W.~Fulton, \emph{Intersection theory}, Ergebnisse der Mathematik und ihrer
  Grenzgebiete (3), vol.~2, Springer-Verlag, Berlin, 1984.

\bibitem{fultonmacpherson:categorical}
W.~Fulton and R.~MacPherson, \emph{Categorical framework for the study of
  singular spaces}, Mem. Amer. Math. Soc. \textbf{31} (1981), no.~243.

\bibitem{goreskymacpherson:intersectionII}
M.~Goresky and R.~MacPherson, \emph{Intersection homology. {II}}, Invent. Math.
  \textbf{72} (1983), no.~1, 77--129.

\bibitem{hiller:geometry}
H.~Hiller, \emph{Geometry of {C}oxeter groups}, Research Notes in Mathematics,
  vol.~54, Pitman (Advanced Publishing Program), Boston, Mass., 1982.

\bibitem{hotta:local}
R.~Hotta, \emph{A local formula for {S}pringer's representation}, Algebraic
  groups and related topics ({K}yoto/{N}agoya, 1983), Adv. {S}tud. {P}ure
  {M}ath., North-Holland, Amsterdam, 1985, pp.~127--138.

\bibitem{kashiwarashapira:sheaves}
M.~Kashiwara and P.~Shapira, \emph{Sheaves on manifolds}, Springer-Verlag,
  1990.

\bibitem{kazhdanlusztig:topological}
D.~Kazhdan and G.~Lusztig, \emph{A topological approach to {S}pringer's
  representations}, Adv. in Math. \textbf{38} (1980), 222--228.

\bibitem{maclane:categories}
S.~Mac Lane, \emph{Categories for the working mathematician}, Graduate Texts in
  Mathematics, vol.~5, Springer-Verlag, New York-Berlin, 1971.

\bibitem{lusztig:green}
G.~Lusztig, \emph{Green polynomials and singularities of unipotent classes},
  Adv. in Math. \textbf{42} (1981), 169--178.

\bibitem{rossmann:picard}
W.~Rossmann, \emph{Picard-{L}efschetz theory for the coadjoint quotient of a
  semisimple {L}ie algebra}, Invent. Math. \textbf{121} (1995), no.~3,
  531--578.

\bibitem{spaltenstein:reflection}
N.~Spaltenstein, \emph{On the reflection representation in {S}pringer's
  theory}, Comment. Math. Helv. \textbf{66} (1991), no.~4, 618--636.

\bibitem{springer:trig}
T.A. Springer, \emph{Trigonometric sums, {G}reen functions of finite groups and
  representations of {W}eyl groups}, Invent. Math. \textbf{36} (1976),
  173--207.

\bibitem{steinberg:desingularization}
R.~Steinberg, \emph{On the desingularization of the unipotent variety}, Invent.
  Math. \textbf{36} (1976), 209--224.

\bibitem{tanisaki:twisted}
T.~Tanisaki, \emph{Twisted differential operators and affine {W}eyl groups}, J.
  Fac. Sci. Univ. Tokyo Sect. IA Math. \textbf{34} (1987), no.~2, 203--221.

\end{thebibliography}
\bibliographystyle{amsplain}

\providecommand{\bysame}{\leavevmode\hbox to3em{\hrulefill}\thinspace}
\providecommand{\MR}{\relax\ifhmode\unskip\space\fi MR }
\providecommand{\MRhref}[2]{%
  \href{http://www.ams.org/mathscinet-getitem?mr=#1}{#2}
}
\providecommand{\href}[2]{#2}


\end{document}